
\documentstyle{amsppt}
\input epsf
\def\cline#1{\leftline{\hfill#1\hfill}}
\NoRunningHeads
\hsize=16cm
\vsize=21cm
\magnification=\magstephalf
\abovedisplayskip=4pt
\belowdisplayskip=4pt
\document
\define\({\left(}
\define\){\right)}
\define\[{\left[}
\define\]{\right]}

\define\eps{\varepsilon}

\define\imbright#1{\overset{\dsize{#1}}\to{\hookrightarrow}}
\define\imbleft#1{\overset{\dsize{#1}}\to{\hookleftarrow}}

\define\al{\alpha}

\define\si{\sigma}
\define\Si{\Sigma}
\define\ga{\gamma}
\define\Ga{\Gamma}

\define\de{\delta}
\define\De{\Delta}

\define\nin{\noindent}
\define\ninpar{\par \noindent}

\define\nlind{\nl \indent}
\define\nl{\newline}
\define\what{\widehat}

\define\wtl{\widetilde}
\define\binv#1{^{(#1)}}
\define\od#1{\frac d{d#1}}
\define\pd#1{\frac \partial{\partial#1}}
\define\PD#1#2{\frac {\partial #1}{\partial#2}}
\define\torus#1{\bold C^{*{#1}}}
\define\ctor#1#2{\bold C_{#1}^{*{#2}}}

\define\weight{\operatorname{weight}\/}

\define\val{\operatorname{val}}

\define\Cone{\operatorname{Cone}\/}

\define\Vol{\operatorname{Vol}\/}

\define\degree{\operatorname{degree}\/}

\define\wrt{with respect to~}
\define\Ts{Tschirnhausen}
\define\TS{Tschirnhausen approximate}

\define\QED{~~Q.E.D.}

\define\bC{\bold C}
\define\bN{\bold N}

\define\bP{\bold P}

\define\bZ{\bold Z}
\define\bR{\bold R}

\define\mapright#1{\smash{\mathop{\longrightarrow}\limits^{\dsize{#1}}}}
\define\mapleft#1{\smash{\mathop{\longleftarrow}\limits^{\dsize{#1}}}}
\define\mapdown#1{\Big\downarrow\rlap{$\vcenter{\hbox{$#1$}}$}}
\define\Trunc{\operatorname{Trunc}\/}
\define\mapdownn#1#2{\llap{$\vcenter{\hbox{$#1$}}$}\Big\downarrow\rlap{$\vcenter{\hbox{$#2$}}$}}
\define\mapup#1{\Big\uparrow\rlap{$\vcenter{\hbox{$#1$}}$}}
\define\mapupp#1#2{\llap{$\vcenter{\hbox{$#1$}}$}\Big\uparrow\rlap{$\vcenter{\hbox{$#2$}}$}}
\define\defby{:=}
\define\eqby#1{\overset {#1}\to =}
\define\inv{^{-1}}

\define\ccfrac#1#2#3{#1 -\dps{1\over\dps{#2-_{\ddots_{
\dps-{1\over\dps{#3}}}}}}}
\define\osmile#1#2{\overset{\overset #1 \to \smile} \to #2}
\define\olor#1#2{\overset{\overset #1 \to \lor} \to #2}

\define\usim#1{\underset{#1} \to \sim}

\parskip=6pt
\parindent=8mm
\baselineskip=14 true pt plus 3pt
\topmatter
\title\nofrills{ Geometry of plane curves via Tschirnhausen resolution tower}
\endtitle
\author Norbert A'CAMPO and  Mutsuo OKA \endauthor
\NoRunningHeads
\address\newline Norbert A'Campo: Mathematisches Institut der Universit\"at,
Rheinsprung 21, CH-4051, Basel,\newline
\phantom{aaaaaaaaaaaaaaaa} Switzerland
\newline
Mutsuo Oka:
Department of Mathematics,
Tokyo Institute of Technology,
Oh-Okayama,\newline\phantom{aaaaaaaaaaaaaaaa}
 Meguro-ku, Tokyo 152, Japan
\endaddress
\email Norbert A'Campo:~nacampo\@math.unibas.ch\newline
\phantom{email address address }
Mutsuo Oka: oka\@math.titech.ac.jp\endemail
\thanks{This work was done when the first author was visiting
the Department of  Mathematics of the Tokyo
Institute of Technology in the fall
 of 1993. We thank the Dept. of Math. of T.I.T.
for their support and hospitality.}
\endthanks
\endtopmatter
\nin
{\bf Table of contents}\nl
\S 1. Introduction\nl
\S 2. \TS ~polynomials of a monic polynomial\nl
\S 3. Toric modification and strict transforms\nl
\S 4. \Ts~ resolution tower for an irreducible germ\nl
\S 5. The zeta-function of the monodromy\nl
\S 6. Conditions implying equi-singularity\nl
\S 7. An example of an equi-singular family\nl
\S 8. The equi-singularity at infinity and the Abhyankar-Moh-Suzuki theorem
\ninpar
\S 1. {\bf Introduction.}
The weight vectors  of a resolution tower of
toric modifications for an irreducible
germ of a plane curve $C$
carry enough information to read  off
 invariants such as the  Puiseux pairs, multiplicities, etc [26].
However, each step  of the inductive
construction of  a tower of toric modifications
depends on a choice of the modification local coordinates.
This  ambiguity makes it difficult to
study  the  equi-singularity problem of
 a  family
of germs of plane curves or to study a global curve.
It is the purpose of this paper
to make a canonical  choice of
 the modification local coordinates $(u_i,v_i)$ ( Theorem (4.5)), and to obtain
 a canonical sequence of germs of curves
$\{C_i; i=1,\dots, k\}~(C_k=C)$ such that the local knot of the curve $C_i$ is
a compound torus knot around  the local knot of the curve $C_{i-1}$. We will
show that the local equations $h_i(x,y)$ of the the germs $\{C_i; i=1,\dots,
k\}~(C_k=C)$ are the
\TS ~polynomials of the local equation $f(x,y)$ for $C$, provided  that
$f(x,y)$ is a  monic polynomial in $y$.
\nlind
The importance of the \TS~ polynomials  was first observed by
Abhyan-\nl kar-Moh [2,3], and our work is very much influenced by them.
However our result gives not only a geometric interpretation of [2,3] but
also  a new method to study the equi-singularity problem for a given
family of germs of irreducible plane curves $f(x,y,t)=0$
whose \TS~ polynomials $h_i(x,y),~i=1,\dots, k-1$
do not depend on $t$.
\nlind
In \S 6, we show that a family of germs of plane curves $\{f_t(x,y)=0\}$
with  \TS~polynomials $h_i(x,y),~i=1,\dots, k-1$ not depending upon $t$ and
satisfying an additional intersection condition is  equi-singular (Theorem
(6.2)).
 In \S  8, we will give a new proof and
a generalization of the Abhyankar-Moh-Suzuki theorem
from the viewpoint of the equi-singularity at infinity
(Theorems (8.2),(8.3),(8.7)).
\par
\nin
\S 2. {\bf  \TS~ polynomials of a monic polynomial.}
Let $f(y)=y^n+\sum_{i=1}^n c_iy^{n-i}$
be a monic polynomial
in $y$ of degree $n$ with coefficients in
 an integral domain $R$ which contains the field of rational numbers $\bold Q,$
 and let $a$ be a positive integer such that
$a$ divides $n.$  The
{\it   $n/a$-th \TS~polynomial}
(or {\it the $n/a$-th \TS~ root})
 of $f(y)$ is the monic polynomial
$h(y)\in R[y]$ of degree $a$ such that
$\degree(f(y)-h(y)^{n/a})<n-a$.
The coefficients of
$h(y)=y^{a}+\sum_{i=1}^a \al_iy^{a-i}$
are  inductively determined by:
$\al_0=1$ and
$c_i(x)=\sum_{j_1+\cdots+j_{\ell }=i} \al_{j_1}(x)\cdots\al_{j_{\ell }}(x)$
for $i=1,\dots, a$.
The coefficient
$\al_j$ is a weighted homogeneous polynomial of degree
$j$ in the variables
$ c_{1},\cdots, c_{a}$
with $\weight(c_{j})=j,\quad 1\le j\le a$. In our application $R$ will be the
ring $\bC\{x\}$ or $\bC[x].$
For further detail, we refer to [2,3].
{}From  the  Euclidean division algorithm it follows
\proclaim{Proposition (2.1)} Let $h(y)\in R[y]$ be monic of
degree $a$ in $y,$ and let $P(y) \in  R[y]$ such that
$sa\le \deg_y P(y) <(s+1) a$.
Then there exits a unique expansion, called the Euclidian  expansion,
$P(y)=\sum_{i=0}^s \al_i(y) h(y)^{s-i},$
where
$\al_i(y)\in  R[y], i=0,\dots, s$ satisfy
 $\deg_y \al_i(y)< a$.
In  particular, we can expand $f(y)$ with respect to its $n/a$-th \TS
{}~polynomial, as
$f(y)=h(y)^{n/a}+\sum_{i=2}^{n/a}c_i(y) h(y)^{n/a-i}$,
$ \deg_y c_i(y)<a.$
If $f(y)=f(x,y) \in \bC[x][y]$,  the coefficients
$c_i(y)=c_i(x,y)$ are also polynomials in $x$ and $y$.
\endproclaim
The second assertion is immediate from the Euclidian expansion of
$f-h^{n/a}$.
We call the above expansion {\it the $n/a $-th ~ Tschirnhausen expansion} of
$f(x,y)$. The  expansion of $P(x,y)$ \wrt~ $h(x,y)$ will also be called {\it
Tschirnhausen expansion} if $h(x,y)$ is a \TS ~polynomial.
 \TS~ polynomials behave hereditary in the following sense.
\proclaim{Proposition (2.2)}
Assume that $a,b\ge 2$ are  integers such that
$ ab|n$. Let $h$ and $h'$ be the $n/a$-th
and $n/ab$-th \TS~ polynomials of $f$ respectively and
let
${h'}=h^{b}+\sum_{i=1}^{b}c_{i} h^{b-i},~\deg_y c_{i}<a$
be the \Ts~ expansion of $h'$ with respect to $h$.
The first coefficient $c_1$  is zero
and
 $h$ is the $ab/a$-th \TS~
polynomial of $h'$.
\endproclaim
\demo{Proof} With  $m:=n/ab$,  we have
$ \deg_y(f-h^{mb})<n-a$ and  $\deg_y(f-{h'}^{m})<n-ab$.
Using the expansion of $h'$ with respect to $h$:
${h'}=h^{b}+\sum_{i=1}^{b}c_{i} h^{b-i},~\deg_y c_{i}<a$, we get
\nl
$${h'}^{m}=\((h^{b}+ c_{1}h^{b-1})+\sum_{i=2}^{b}c_{i} h^{b-i}\)^m
=h^{mb}+ mc_{1}h^{mb-1}+R_1+R_2,
$$
where
$  R_1:=\sum_{i=2}^{m}\binom mi c_1^ih^{mb-i}$ and
$R_2:=\sum_{i=1}^m \binom mi (h^{b}+ c_{1}h^{b-1})^{m-i}\(\sum_{i=2}^{b}c_{i}
h^{b-i}\)^i$.\nl
If $c_1 \not=0$, we would first conclude
$\deg_y R_1<n-a+\deg_y c_1,\quad \deg_yR_2 <n-a,$
and then
$$\align
n-a>\deg_y(f-h^{mb})&=\deg_y(f-{h'}^{m}+mc_{1} h^{mb-1} +R_1+R_2)\\
&=\deg_y(c_{1} h^{mb-1})\ge(mb-1)a= n-a.
\endalign$$
So $c_1=0$ and
it follows
$\deg_y(h'-h^b)=\deg_y\(\sum_{i=2}^{b}c_{i} h^{b-i}\)<ab-a.$
By the uniqueness of the \TS~ polynomial,
the above inequality implies that
$h$ is the $ab/a$-th \TS~polynomial of $h'$.~\QED
\enddemo
The generalized binomial formula:
$(1+z)^r=\sum_{j=0}^\infty \binom{r}{j}z^j$ for $ r>0,$ with coefficients
$\binom{r}{j}={r(r-1)\cdots(r-j+1)}/{j!},$ converges for $|z|<1$. When $r$ is a
rational number $p/q$, the identity:
$\((1+z)^{p/q}\)^q=(1+z)^p$ gives
a recurrent computation of the coefficients of $(1+z)^{p/q}.$
In particular, with $\Trunc^{(\ell)}(1+z)^{p/q}:=\sum_{j=0}^\ell
\binom{p/q}{j}z^j,$
it follows that
$$\val_z\((1+z)^p-(\Trunc^{(\ell)}(1+z)^{p/q})^q\)>\ell\tag 2.2.1$$
For a  real number $x\in \bR$,  denote by $[x]$
the  largest integer $n$ such that $n\le x$.
\proclaim{Lemma (2.3)}
Assume that $a,b,c,d$ are positive integer such that $\gcd(a,b)=1 $ and that
 $d$ divides $ac$. Let
$F(y,z)=(y^a+z^b)^c$ and
$H(y,z)=y^{ac/d} \Trunc^{([c/d])}(1+z^b/y^a)^{c/d}.$
Then  $H$ is the $d$-th \TS~
polynomial of $F(y,z)$ as a polynomial of $y$.
\endproclaim
\demo{Proof}The polynomials $F(y,z)$ and $H(y,z)$ are
weighted homogeneous  of degree $abc$  and $abc/d$
respectively with respect to the weight vector $P={}^t(b,a)$.
In particular, the monomials in $F(y,z)$ and $H(y,z)^d$
have the form $y^{ai}z^{bj}$ with $i+j=c$.
Note also that $\deg_y F(y,z)=ac$,  $\deg_y H=ac/d$ and
$\deg_y(F-H^a)<ac-a[c/d]$ by (2.2.1).
As $ac-ac/d>ac-a[c/d]-a$, this implies the inequality:
$\deg_y\( F(y,z)-H(y,z)^d\)< ac-ac/d.$~\QED
\enddemo
\par\nin
\S 3. {\bf Toric modifications and strict transforms.}\nl
(3.1) {\bf Basic properties of toric modifications} (see [24,27,28,31]).
 Let $(x,y)$ be a fixed system of
local (or global) coordinates of $\bC^2$ at the origin.
Let $N$ be the lattice of integral weights for the monomials in $(x,y)$.  The
weights $E_1(x^ay^b)=a$ and $E_2(x^ay^b)=b$ span the lattice $N,$
and a weight $\al_iE_1+\beta_iE_2$ will be denoted by
the integral column vector ${}^t(\al_i,\beta_i)$.
Let   $N^+$ be the space  of positive weight  vectors of $N,$ and similarly let
$N_{\bR}^+$ be the positive cone in $N_{\bR}:=N\otimes_{\bZ} \bR$.
A simplicial cone  subdivision $\Si^*$ of $N_{\bR}^+$ is a sequence
$\(T_1,\dots, T_m\)$ of primitive  weights in $N^+,$ called the vertices, such
that $T_0=E_1,T_{m+1}=E_2$ and
$\det(T_i,T_{i+1})=\det_{\{E_1,E_2\}}(T_i,T_{i+1}) \ge 1$ holds.
 The $m+1$  cones
$\Cone(T_i,T_{i+1}):=\{ tT_i+s T_{i+1}; t,s \ge 0\}$,
 $i=0,\dots, m,$ cover without overlap the  cone  $N_{\bR}^+$.
The subdivision $\Si^*$ is called {\it  regular }
if
$\det(T_i,T_{i+1})=1$ for  each $i=0,\dots, m$.
Let $\si_i$ be the integral matrix mapping $E_1$ to $T_i$ and $E_2$ to
$T_{i+1}$.
\nlind
 Using a birational mapping
$\phi_M:\bC^2 \to \bC^2$,
$\phi_M(x,y)=(x^ay^b,x^cy^d)$ for an  integral unimodular matrix
$M=\pmatrix a&b\\c&d\endpmatrix$,
the  toric modification $p: X\to \bC^2$ associated with a regular simplicial
cone subdivision $\Si^*$
 is defined as follows.
The non-singular complex manifold $X$ is covered with $m+1$ so called {\it
toric coordinate charts}
$\{\bC_{\si_{i}}^2,(x_{\si_{i}},y_{\si_{i}})\},~i=0,\dots, m,$ where  points
$(x_{\si_i},y_{\si_i})\in \bC_{\si_i}^2$ and
$(x_{\si_j},y_{\si_j})\in \bC_{\si_j}^2$ are identified
if and only if the birational map
$\phi_{\si_{j}\inv \si_{i}}$ is defined at the point
$(x_{\si_{i}},y_{\si_{i}}) \in \bC_{\si_i}^2$ and
$\phi_{\si_{j}\inv \si_{i}}(x_{\si_i},y_{\si_i})=(x_{\si_{j}},y_{\si_{j}})$.
The  morphism
$\pi_{\si_i}: \bC_{\si_i}^2\to  \bC^2$  defined by
$
\pi_{\si_i}(x_{\si_{i}},y_{\si_{i}})=\phi_{\si_i}(x_{\si_{i}},y_{\si_{i}})$ are
compatible with the identifications and define a proper birational analytic map
$p: X\to \bC^2$ (see [31,28]).
A toric modification is a composition of finite blowing-ups
(see [17]).
 The exceptional divisor
 $p\inv(O)$ is the union of  $m$ rational curves
$\{\what E(T_i);i=1,\dots,m\}$ and
 each one is  covered by {\it its left chart}
$\bC_{\si_{i-1}}^2$ and {\it its right chart} $\bC_{\si_i}^2$ and defined by
the equations
$\{x_{\si_i}=0\}$, $\{y_{\si_{i-1}}=0\}$. Thus only $\what E(T_i)$ and $\what
E(T_{i+ 1})$
 intersect transversely at the origin of the chart $\bC_{\si_i}^2$.
The non-compact divisor
$\what E(E_1):=\{x_{\si_{0}}=0\}$ and
$\what E(E_2):=\{y_{\si_m}=0\}$  map isomorphically
onto
the axis $x=0$ and $y=0$.
\par\nin
(3.2) {\bf Admissible toric modifications.}
Let  $f(x,y)=\sum a_{\alpha,\beta}x^{\alpha}y^{\beta}$ be the Taylor
expansion  of a germ of a holomorphic function f
with $f(O)=0$. The Newton polygon $\Gamma_+(f;(x,y))$ of $f(x,y)$ is the
 convex hull in  $N_{\bR}^*$ of
$\{(\alpha+s,\beta+t) \in \bR^2 ;  a_{\alpha,\beta}\not=0,s \ge 0,t \ge 0\}$
and  the Newton boundary
$\Ga(f;(x,y))$ is the union of the compact faces of $\Gamma_+(f;(x,y))$    (see
[24,25,27] for instance).
The Newton boundary $\Ga(f;(x,y))$ contains only a finite number of faces
of dimension one.
Each positive weight vector $P={}^t(p,q)\in N^+$
defines a non negative function  on $\Gamma_+(f;(x,y))$, for which we denote by
$d(P;f)$ its minimal value  and by
$\De(P;f)$  the face or the vertex where this minimal value is taken.
We consider on  $N^+$ the equivalence relation:
$P\sim Q$ iff $\De(P;f)=\De(Q;f)$. {\it The dual Newton diagram}
$\Ga^*(f;(x,y))$ of $f(x,y)$
is  the conical subdivision  of $N^+$ given by the equivalence classes. Let
$P_i={}^t(a_i,b_i) \in N^+,i=1, \dots ,m$ be the ordered list of   primitive
weight vectors such that
 $\De(P_i;f)$ is  the list of the one-dimensional faces of $\Ga^*(f)$ and
$\det(P_i,P_{i+1})= a_ib_{i+1}-a_{i+1}b_i>0$,  $i=1,\dots, m-1$.
The {\it face function}
$f_{P_i}(x,y):=\sum_{(\al,\beta)\in \De(P_i;f)} a_{\al,\beta} x^\al y^\beta$
 admits a product decomposition
$f_{P_i}(x,y)=c_ix^{r_i}y^{s_i}
\prod_{j=1}^{k_i}(y^{a_i}-\gamma_{i,j}x^{b_i})^{\nu_{i,j}}$
with distinct non-zero complex numbers
$\gamma_{i,1},\ldots,\gamma_{i,k_i}$.
Recall that $f(x,y)$ is {\it non-degenerate} if and only if
$\nu_{i,j}=1$ for any $i,j$.
The partial sum
$\Cal N(f)(x,y)=\sum '' a_{\al,\beta} x^\al y^\beta$
over  all $(\al,\beta) \in \Ga(f;(x,y))$ is the {\it  Newton principal part}
$\Cal N(f)(x,y)$.
\nlind
 A regular simplicial cone subdivision
$\Si^*$ with vertices $\{T_0=E_1,T_1,\dots, T_\ell, T_{\ell+1}=E_2\}$
 is called {\it admissible} for $f(x,y)$
 if $\Si^*$ is a refinement
of the dual Newton diagram $\Ga^*(f;(x,y))$, meaning $P_i={}^t(a_i,b_i) \in
\{T_0,T_1,\dots, T_\ell, T_{\ell+1}\},i=1, \dots ,m$.
Note that $\Si^*$ is admissible for $f(x,y)$ if and only if
$\De(T_j;f)\cap \De(T_{j+1};f)\ne \emptyset,j=0,\dots, \ell$.
 The corresponding toric modification
$p : X\to \bC^2$ is called {\it  admissible } for $f(x,y)$.\nlind
Basic properties of  admissible
toric modifications are:\nl
(3.2.A) The
divisor $\what E(T_{j})$  meets the proper transform
$\wtl C$ if and only if $T_j$ is a primitive weight $P_i$.
\nl
(3.2.B) The  divisor $\what E(P_i)$ intersects $\wtl C$ at
$k_i$ points. In the right toric  chart
$\{\bC_{\si_j},(x_{\si_j},y_{\si_j})\},\si_j=\Cone (T_{j},T_{j+1}),P_i=T_j$
of $\what E(P_i)$,
the intersection   $\wtl C\cap \what E(P_i)$ is $\{(0,\ga_{i,1}),\dots,
(0,\ga_{i,k_i})\}$.\nl
(3.2.C) The divisor of the pull back  $p^*f$
of the function  $f$ is given by
$$(p^*f)=\sum_{i=1}^m\sum_{\ell=1}^{k_i}\wtl C_{i,\ell}+
\sum_{j=0}^{\ell+1}d(T_{j};f)\what E(T_{j})
$$
where  $\wtl C_{i,\ell}$ is the union of components of $\wtl C$ which
 pass through $(0,\ga_{i,\ell})$. \nl
(3.2.D)  If $f(x,y)$ is irreducible as
a germ of a function at the origin, then
$m=1$ and $k_1=1$.\nl
(3.2.E) If $f$ is non-degenerate,
the curve $\wtl C_{i,j}$ is  smooth and
$\wtl C_{i,j}$ intersects transversely with
 $\what E(P_i)$.
Thus, if $f(x,y)$ is non-degenerate, the modification $p$ is a
good resolution of $f(x,y)$ (see [17]).
\ninpar
(3.3) {\bf Intersection multiplicity with a reduced irreducible germ.}
Let $C=\{f(x,y)=0\}$  be a  reduced irreducible germ of a   curve. The defining
function admits for a weight $P_1={}^t(a_1,b_1)$ an initial expansion
$f(x,y)=(y^{a_1}+\xi_1 x^{b_1})^{A_2}+\text{(higher terms)}$
with $ \xi_1\not=0$ and $\gcd(a_1,b_1)=1,$
 where ``higher terms'' collects the monomials of $P_1$-degree strictly
greater than $a_1b_1A_2.$
Let $C'$ be another (not necessarily irreducible)
germ of a curve
defined by $C'=\{(x,y)\in U;g(x,y)=0\}$.
Let $p:X\to \bC^2$ be a toric modification admissible both for $C$ and $C'$,
and let $\Xi_1$ be the intersection point of $\wtl C$ and $E(P_1)$.
Then
\proclaim{Proposition (3.3.1)(Lemma (7.12),[27])}
The intersection multiplicity of $C$ and $C'$
at the origin is
$I(C,C';O)=d(P_1;g)A_2 + I(\wtl C,\wtl C';\Xi_1)$.
The  term $I(\wtl C,\wtl C';\Xi_1)$ vanishes if and only if
$g_{P_1}(x,y)$ is not divisible by $(y^{a_1}+\xi_1 x^{b_1})$.
If $g(x,y)$ has for a primitive weight vector
$P_1'={}^t(a_1',b_1')$ the initial expansion
$g(x,y)=(y^{a_1'}+\xi_1' x^{b_1'})^{A_2'}+\text{(higher terms)}$,
then $d(P_1;g)A_2=\min(a_1b_1',a_1'b_1)\times A_2A_2'$
and moreover
$I(\wtl C,\wtl C';\Xi_1)=0$
 if and only if either $P_1\ne P_1'$
or $P_1=P_1'$ and  $\xi_1\ne \xi_1'$. \endproclaim
\comment\demo{Proof}
We can write
$$\align
(p^* f)&=\wtl C+a_1b_1A_2 \what E(P_1)+D\\
(p^*g)&=\wtl C'+d(P_1;g)\what E(P_1)+D'\endalign$$
where $D,D'$ are linear combinations of the compact divisors other than
$\what E(P_1)$. Note that \nl
$I(C,C';O)=(p^* f)\cdot (p^*g)$ where the right side is the intersection
number in a small neighbourhood of the exceptional divisor $p\inv (O)$.
Then
the assertion follows from the equalities $\wtl C\cdot \what E(P_1)=A_2$
and
$(a_1b_1A_2 \what E(P_1)+D)\cdot (p^*g)=0$.\enddemo
\ninpar
\endcomment
\nin
(3.4) {\bf A resolution tower of toric modification for
an irreducible germ.}
Let $C$ be an irreducible germ of a plane curve   and let
$\Cal T =
\{ X_k\mapright{p_k}X_{k-1}\mapright{p_{k-1}}\dots
\to X_1\mapright{p_1}X_0\}$
be a sequence of non-trivial toric modifications
where
each $p_{i+1}:X_{i+1}\to X_{i}$ is the toric modification
associated with a regular  simplicial cone subdivision $\Si_i^*$ of the cone
$N_{\bR}^+$ in the space of weights for
a local system  of coordinates $(u_i,v_i)$ of $X_{i}$, centered at
 the center $\Xi_{i}\in X_i$ of the modification $p_{i+1}$.
Let $E_{i,1},\dots, E_{i,s_i}$ be the exceptional divisors of
$p_i:X_i\to X_{i-1}$.   Abusely,
we denote by the same $E_{i,j}$,  the strict transforms of $E_{i,j}$ to
$X_\ell$ for any $\ell\ge i$.
Thus the exceptional divisors of  the modification
$\Phi_k\defby p_1\circ\dots\circ p_k:X_k\to X_0$
are $\{E_{i,j}\},1\le i\le k, 1\le j\le s_i$.
Denote by $\Xi_{i}\in E_{i,\beta_i}\cap C^{(i)}$
the preimage of the singularity in the strict transform $C^{(i)}$ of $C$ to
$X_i$.
 We call $\Cal T$  {\it a resolution tower of  admissible toric modifications}
if the following conditions are satisfied ([27]).\nl
(i) $X_0$ is an open neighborhood of the origin $O$ of $\bC^2$,
$(u_0,v_0)=(x,y)$ and $\Xi_0=O$.
\nl
(ii) The modification $p_{i+1}:X_{i+1}\to X_{i}$ is  non-trivial
 and  admissible  for $\Phi_i^* f(u_i,v_i),\quad i>0$.
\nl
(iii) The coordinate  $u_i$ is simply the restriction
$u_i=x_{\si_i}|{W_i}$ of the coordinate $x_{\si_i}$ of the right toric chart of
$E_{i,\beta_i}$ to  a neighborhood $W_i$ of $\Xi_i$.
\nl
(iv) $p_i(\Xi_i)=\Xi_{i-1}$.
\nl
(v) The composition
 $\Phi_k: X_k\to  X_0$ is a good resolution of $C$.
\nlind
The weight vectors
$P_i={}^t(a_i,b_i)$  corresponding
to the exceptional divisors $E_{i,\beta_i}$ for $i=1,\dots, k$
 are
 {\it the weight vectors}  {\it of the tower} ([27]).
 If the tower $\Cal T$ is  admissible for $C$, there exists for $i=0,\dots,
k-1$ non-zero complex numbers $\xi_i\in \bC$
so that
$
C^{(i)}=\{(u_i,v_i)\in W_i;
{}~(v_i^{a_{i+1}}+\xi_{i+1}u_{i}^{b_{i+1}})^{A_{i+2}}+
\text{(higher terms)}=0\}$
 where $C^{(0)}=C$ and $A_j=a_j\cdots a_k,j\le k$ and $A_{k+1}=1$.
\par
Let $D=\{(x,y)\in \bC^2; g(x,y)=0\}$
be an irreducible, not necessarily reduced,
germ of a plane curve at the origin of $\bC^2=X_0$ and let $D^{(i)}$ be the
strict transform of $D$ to $X_i$. If
$D$ has  the same toric tangential direction of depth $\theta$
with $C$ \wrt~ $\Cal T,$ i.e. if $\Xi_i\in D^{(i)}$ for $i\le \theta$
and $\Xi_{\theta+1}\notin D^{(\theta+1)},$  there exist a
non-zero complex number $\xi_{\theta+1}'$, a positive
integer ${A_{\theta+2}}'$ and a primitive weight
vector $P_{\theta+1}'\defby {}^t(a_{\theta+1}',b_{\theta+1}')$
such that
$$
D^{(i)}=\cases
\{(u_i,v_i)\in W_i;
{}~(v_i^{a_{i+1}}+\xi_{i+1}u_{i}^{b_{i+1}})^{A_{i+2}'}+
\text{(higher terms)}=0\},\quad &i <\theta\\
\{(u_\theta,v_\theta)\in W_\theta;
{}~(v_\theta ^{a_{\theta +1}'}+\xi_{\theta +1}'
u_{\theta }^{b_{\theta +1}'})^{A_{\theta  +2}'}+
\text{(higher terms)}=0\},\quad &i =\theta\endcases\tag3.4.1
$$
where
$ A_j'=a_j\cdots a_\theta a_{\theta+1}'A_{\theta+2}'$, $j\le \theta+1$.
If $P_{\theta+1}'={}^t(1,0)$, the transform $D^{(\theta)}$ is defined
by $\{v_{\theta+1}^{A_{\theta +2}'}=0\}$ since $D$ is  irreducible.
The case $P_{\theta+1}'={}^t(0,1)$ does not occur as $\{u_{\theta+1}=0\}$
is nothing but $\what E(P_{\theta})$.
Put
$$I(P_{\theta+1},P_{\theta+1}')
:=\cases
\min(a_{\theta+1}b_{\theta+1}',a_{\theta+1}'b_{\theta+1}),
\quad&\text{if~} a_{\theta+1}b_{\theta+1}',a_{\theta+1}'b_{\theta+1}>0\\
b_{\theta+1},\quad &\text{if~}
P_{\theta+1}'={}^t(1,0)\endcases$$
 By  induction, using
 Proposition (3.3.1), we get
\proclaim{Lemma (3.4.2) ([27])} Assume that
$D$ has  the same toric tangential direction of depth $\theta$
with $C$ \wrt~ $\Cal T$. Under the  assumption  (3.4.1),
 the local intersection multiplicity is
$$
I(C,D;O)=\sum_{i=1}^\theta a_ib_i A_{i+1}A_{i+1}'
+I(P_{\theta+1},P_{\theta+1}')\times A_{\theta+2}
A_{\theta+2}'$$
\endproclaim
Let $D_1,\dots ,D_r$
be the irreducible components of a reducible  plane curve germ $D$.  We say
that
the reducible germ $D$ has {\it the same toric tangential direction of depth
$\theta$
with $C$ \wrt~ $\Cal T$} if $\Xi_i\in D_j^{(i)}$ for any
$j=1,\dots,r$ and  $i\le \theta$
and $\Xi_{\theta+1}\notin D_{j_0}^{(\theta+1)}$ for some $j_0$.
\par
\nin
\S 4. {\bf A Tschirnhausen resolution  tower for an irreducible
germ.}
\proclaim{Lemma (4.1)} Let $p:X\to \bC^2$ be a  toric modification
with respect to a regular simplicial cone subdivision $\Si^*$
of $N^+$.
 Let $\si=\Cone (P,P')$ be a  cone   in $\Si^*$  and  $g(x,y) \in \bC\{x,y\}$,
such that
$\De(P;g)$ is a vertex. At each point
$\Xi\in \what E(P)-\bigcup_{Q\ne P}\what E(Q)$
the function $p^* g/x_{\si}^{d(P;g)}$
is  a unit.
\endproclaim
\demo{Proof}Let $\{(\nu_1,\nu_2)\}=\De(P;g)$
and  $c\ne 0$ be the coefficient of $x^{\nu_1}y^{\nu_2}$
in  $g(x,y)$.
Then
$p^* g=x_{\si}^{d(P;g)}y_{\si}^{d(P';g)}
\{c  y_\si^\al +x_\si  g'(x_{\si},y_{\si})\}$
for some analytic function $g'(x_{\si},y_{\si})$ and $\al\ge 0$. Moreover,
$\al=0$ iff $\De(P';g)\supset \De(P;g)$. In conclusion, $p^*
g/x_{\si}^{d(P;g)}$
is a unit at $\Xi$ since $y_\si$ is.
\QED\enddemo
In particular, if $\Ga(g;(x,y))\subset \{(\nu_1,\nu_2);\nu_2<a\}$ or  if
$g(x,y) \in \bC\{x\}[y]$ and
$\deg_y g< a$, the face $\De(P;g)$ is a vertex, and the lemma applies.
\ninpar
{\bf A. \Ts ~ resolution Tower.}\nl
{\bf (4.2)} Let $f(x,y) \in \bC\{x\}[y]$  be monic of
 degree $n$ and irreducible   with the initial expansion
$f(x,y)=(y^{a_1}+\xi_1 x^{b_1})^{A_2}+\text{(higher terms)}$,
$a_1>1$,
for the primitive weight vector $P_1={}^t(a_1,b_1)$ with $n=a_1A_2$.
The $n/a$-th
\TS~ polynomial $H_a(x,y)$ is a monic polynomial of degree $a$ in $y$ and
defines  at the origin the germ of the curve
$D_a:=\{H_a(x,y)=0\}.$
\par\noindent
{\bf  (4.3) First Observation.}
Let $p_1:X_1\to \bC^2$ be an
admissible toric modification with respect to a regular simplicial cone
subdivision
$\Si_0^*$ for  $f(x,y)$.
The strict transform
$C^{(1)}$ of $C$ to $X_1$ intersects  only with $\what E(P_1)$, say at the
point
 $\Xi_1.$  In the
chart $\bC_{\si_1}^2,$ where $\si_1=(P_1,P_1')$, $P_1'={}^t(a_1',b_1')$ is  the
right cone
of  $\what E(P_1),$ we have
$\Xi_1=(0,-\xi_1)$.
Put
$h_1(x,y):=H_{a_1},~ C_1:=D_{a_1}.$ For a multiple $a$ of $a_1$ with $a|n,$
the
$A_2$-th (resp.  $n/a$-th )
\TS~ polynomial of $(y^{a_1}+\xi_1 x^{b_1})^{A_2}$
 is the face function $h_{1 P_1}$ (respectively $H_{a P_1}$), hence
$h_1$ and $H_a$ can be written as:
$$\cases
&h_1(x,y)=y^{a_1}+\xi_1 x^{b_1}+\text{(higher terms)}\\
&H_a(x,y)=(y^{a_1}+\xi_1 x^{b_1})^{a/a_1}+\text{(higher terms)},\quad
\text{if}\quad a_1|a
\endcases
\tag 4.3.1
$$
In particular, $h_1(x,y)$ is non-degenerate.
As
$p_1^*(y^{a_1}+\xi_1x^{b_1})=x_{\si_1}^{a_1b_1}y_{\si_1}^{a_1'b_1}(y_{\si_1}+\xi_1)$,
we can write
$
p_1^*h_1(x_{\si_1},y_{\si_1})=
x_{\si_1}^{a_1b_1}y_{\si_1}^{a_1'b_1}\((y_{\si_1}+\xi_1)+
x_{\si_1}R(x_{\si_1},y_{\si_1})\),R(x_{\si_1},y_{\si_1}) \in
C\{x_{\si_1},y_{\si_1}\}.$
The  functions
$u_1=x_{\si_1}$, $ v_1=p_1^*h_1/x_{\si_1}^{a_1b_1}
=y_{\si_1}^{a_1'b_1}\((y_{\si_1}+\xi_1)+
x_{\si_1}R(x_{\si_1},y_{\si_1})\)
$
give a system of  coordinates $(u_1,v_1)$  in a neighbourhood $W_1$ of
  $\Xi_1.$
The
strict transform $C_1^{(1)}$ of $C_1$   to $X_1$
intersects only with $\what E(P_1)$ and
$p_1^*h_1=u_1^{a_1b_1} v_1,~
C_1^{(1)}=\{v_1=0\},$
so $C_1$ is irreducible and $p_1$ is a good resolution of $C_1$.
\nl
If  $A_2=1$, we have $f=h_1$
and we have nothing  to do further.
If  $A_2\ge 2,$
 the pull back $p_1^*f(u_1,v_1)$ has an initial expansion
$$
p_1^*f(u_1,v_1)=u_1^{m_1(f)}(v_1^{a_2}+\xi_2 u_1^{b_2})^{A_3}+
\quad\text{(higher terms)}\tag 4.3.2$$
with primitive weight vector $P_2={}^t(a_2,b_2)$
where   the multiplicity  of
$\Phi_1^*f$ on $E_1$ is  $m_1(f)=a_1b_1A_2$ by
 (4.2.1). Note also $A_2=a_2A_3$ and
$I(C_1,C;O)=a_1b_1A_2+ b_2A_3$ by   Lemma (3.4.2).
The advantage of the ``\Ts~ coordinates'' is the  inequality $a_2\ge 2$.
In fact, in the  \Ts~ expansion
$f(x,y)=h_1(x,y)^{A_2}+\sum_{j=2}^{A_2}c_{j}(x,y)h_1(x,y)^{A_2-j}$ of $f(x,y)$
\wrt~ $h_1$
we have $ c_{j}(x,y)\in \bC\{x\}[y]$ and $\deg_y c_{j}(x,y)<a_1,j=2,\dots,
A_2,$
 so the face $\De(P_1, c_{j})$ is necessarily
a vertex.
Therefore by the definition of the coordinate $(u_1,v_1)$
and Lemma (4.1)  in  a smaller  neighbourhood $W_1$ of $\Xi_1$ the pull-backs
are:
$p_1^*h_1(u_1,v_1)=u_1^{m_1(h_1)}v_1$ with $m_1(h_1)=a_1b_1$ and
$p_1^*c_{j}(u_1,v_1)=u_1^{m_{j}}U_j$
where $ m_{j}=d(P_1,c_{j})$ and
 $U_j$ is a unit for  $j\ge 2$ with $c_j\ne 0$.
If $c_j=0$, we put  $U_j=0$ for simplicity.
Thus  we have
$$p_1^*f(u_1,v_1)=(u_1^{m_1(h_1)}v_1)^{A_2}+
\sum_{j=2}^{A_2}u_1^{m_{j}}U_j(u_1^{m_1(h_1)}v_1)^{A_2-j},\tag 4.3.3$$
hence, with
$Q_0=(m_1(h_1)A_2,A_2)$ and
$ Q_j=(m_{j}+(A_2-j) m_1(h_1),A_2-j),$ the Newton polygon
$\Ga_+(p_1^*f;(u_1,v_1))$
is the convex hull of the sets
$\{Q_0+\bR_+^2\}$ and $\{Q_j+\bR_+^2\},2 \le j \le A_2,c_j\ne 0.$
The Newton principal part $\Cal N(p_1^*f)(u_1,v_1)$
contains
$(m_1(f)+b_2,A_2-a_2)$ by (4.3.2).  It follows, that $(m_1(f)+b_2,A_2-a_2)=Q_j$
for $j=0$ or for some $j \ge 2$, hence
$a_2\ge 2.$ Moreover, if $c_j\ne 0,$ we have $d(P_2; p_1^*c_{j}h_1^{A_2-j})\ge
d(P_2;p_1^*f)$,
with equality
if and only if  $a_2|j$.\nlind
Let $a_1|a$ and $a|n$. The following  \Ts~ expansions start at $j=2$ by
Proposition (2.2):
$$
\cases
H_a&=
h_1^{a/a_1}+\sum_{j=2}^{a/a_1}d_j h_1^{a/a_1-j}\in
\bC\{ x\}[y][h_1],~\deg_y d_j<a_1\\
f&=H_a^{n/a}+\sum_{j=2}^{n/a} c_jH_a^{n/a-j}
\in \bC \{ x\}[y][H_a], ~~\deg_yc_j<a\endcases
\tag4.3.4
$$
By (4.3.3), the principal part of $p_1^*f(u_1,v_1)$
 \wrt~ the weight vector
$P_2$ is
$$p_1^*f_{P_2}(u_1,v_1)=u _1^{m_1(f)}(v_1^{a_2}+\xi_2 u_1^{b_2})^{A_3}
\tag 4.3.5$$
With $R_a:=\sum_{j=2}^{n/a} c_jH_a^{n/a-j}$,
we have that $\deg_y R_a<n-a$ and
therefore the \Ts~expansion of $R_a$
\wrt~ $h_1$  can be written as
$R_a=\sum_{\ell=0}^{A_2-a/a_1-1}\beta_\ell h_1^{\ell}$ for some
$\beta_\ell\in \bC\{x\}[y]$ and $\deg_y\beta_\ell< a_1$.
If $\beta_\ell\ne 0$, we can by Lemma (4.1)
for a  unit $U_\ell$
 and  non-negative integer $\ga_\ell\in \bN$ write $p_1^*(\beta_\ell
h_1^{\ell})=U_\ell u_1^{\ga_\ell}v_1^{\ell}$.
Thus we have for the Newton principal part
$$\deg_{v_1}\Cal N(p_1^*R_a)(u_1,v_1)\le  A_2-a/a_1-1.\tag4.3.6$$
So,  comparing the pull-back
$p_1^*f(u_1,v_1)=p_1^*H_a^{n/a}+\sum_{j=2}^{n/a}U_\ell
 u_1^{\ga_\ell}v_1^{\ell}$ of (4.3.4)
 and (4.3.5),
we see that the monomials
$u_1^{m_1(f)}\times\binom{A_3}{i} (v_1^{a_2})^i(\xi_2 u_1^{b_2})^{A_3-i}$
of $u_1^{m_1(f)}(v_1^{a_2}+\xi_2 u_1^{b_2})^{A_3}$ for
$i>A_3-a/a_1a_2-1/a_2$
 come from $p_1^*H_a^{n/a}$.
The expansions   (4.3.4) and (4.3.1)
give with some analytic functions $g_a,G_a$
$$\cases
&p_1^*H_a(u_1,v_1)=u_1^{m_1(H_a)}(v_1^{a/a_1}+u_1g_a(u_1,v_1))\\
&p_1^*H_a^{n/a}(u_1,v_1)=u_1^{m_1(f)}(v_1^{A_2}+u_1G_a(u_1,v_1))
\endcases\tag4.3.7
$$
 Note that
$m_1(H_a)n/a=m_1(f)$.
Applying the above argument to $p_1^*H_a(u_1,v_1),$  we can conclude that the
Newton boundary $\Ga(p_1^*H_a;(u_1,v_1))$
is situated in the region $\{(\nu_1,\nu_2)\in \bR^2; 0\le \nu_2\le a/a_1\}$
and that $B_a:=(m_1(H_a), a/a_1)$ is the vertex of the left end of
$\Ga(p_1^*H_a;(u_1,v_1))$ by (4.3.7). Note also that $(n/a) B_a=(m_1(f),A_2)$
is the left end of $\Ga(p_1^*f;(u_1,v_1))$ by (4.3.7). Let $\De_a$ be the
first face of $\Ga(p_1^*H_a)$ which contains $B_a$ and let
$Q={}^t(p_2,q_2)$ be the weight vector of $\De_a$.
\proclaim{Assertion (4.3.8)}
The inequality  $q_2/p_2\ge b_2/a_2$ holds.\endproclaim
\demo{Proof}Assuming by contradiction that
$q_2/p_2< b_2/a_2$, we  have
$p_1^*f_{Q}(u_1,v_1)=u^{m_1(f)}v_1^{A_2}$, and we will prove the assertion
by excluding the following three cases:
(a) $d(Q;p_1^*H_a^{n/a})>d(Q;p_1^* R_a)$,
\nl
(b) $d(Q;p_1^*H_a^{n/a})<d(Q;p_1^* R_a)$,
(c) $d(Q;p_1^*H_a^{n/a})=d(Q;p_1^* R_a)$.
Figure (4.3.A) indicates the respective situations.
In case (a), $u^{m_1(f)}v_1^{A_2}=(p_1^* R_a)_Q(u_1,v_1)$ holds,
which is impossible by (4.3.6). The case (b) is impossible
as  $(p_1^*H_a)_Q(u_1,v_1)^{n/a}\ne u^{m_1(f)}v_1^{A_2}$ by the assumption.
If case (c) holds, from (4.3.6) it follows
$(p_1^*H_a)_Q(u_1,v_1)^{n/a}+(p_1^* R_a)_Q(u_1,v_1)\ne 0$, and then
$d(Q;p_1^*H_a^{n/a})=d(Q;p_1^* R_a)=d(Q;p_1^*f)$ and finally the equality
$u^{m_1(f)}v_1^{A_2}=(p_1^*H_a)_Q(u_1,v_1)^{n/a}+(p_1^* R_a)_Q(u_1,v_1).$ \nl
But this  equality is impossible. In fact, let us write
$(p_1^*H_a)_Q(u_1,v_1)=u_1^{m_1(H_a)}v_1^{a/a_1}+
\ga u_1^{\al_1}v_1^{\beta_1}+S(u_1,v_1)$
where $\ga\ne 0$, $0\le \beta_1< a/a_1$ and  $\deg_{v_1}S(u_1,v_1)<\beta_1$ if
$S\ne 0$.
Then
$(p_1^*H_a)_Q(u_1,v_1)^{n/a}=u_1^{m_1(f)}v_1^{A_2}+
n/a\cdot \ga u_1^{\al_1'}v_1^{\beta_1'}+S'(u_1,v_1)$ with
$\deg_{v_1} S'<\beta_1'$
where $\alpha_1'=\al_1+(n/a-1)m_1(H_a)$
and  $\beta_1'=A_2-a/a_1+\beta_1\ge A_2-a/a_1$. On the other hand,
the second term of the right side of the equality has no
monomial $u_1^{\nu_1}v_1^{\nu_2}$ with $\nu_2\ge A_2-a/a_1$.
{}~\QED\enddemo
\midinsert
\cline{\epsffile{key4.3.A.ps}}
\medskip
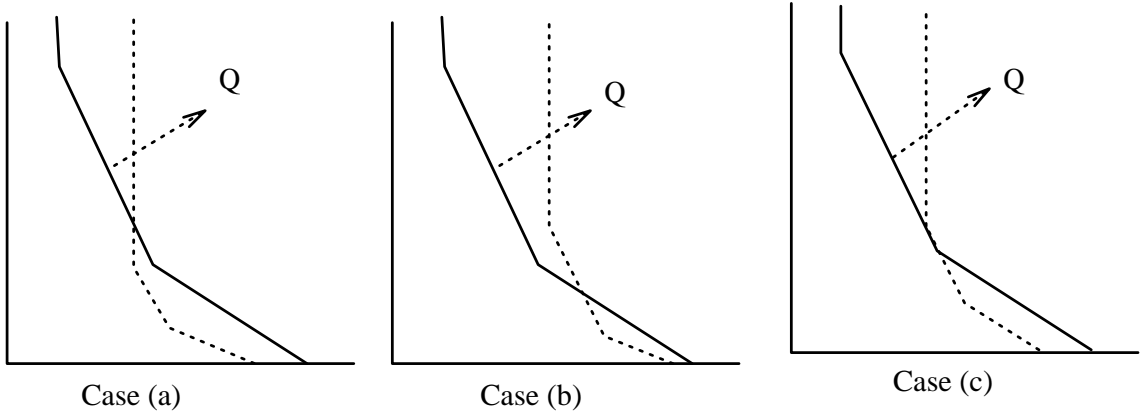
\centerline{Figure (4.3.A)}
\endinsert
See Figure (4.3.A).
Consider the face function $(p_1^*H_a)_{P_2}(u_1,v_1)$ for
the weight vector $P_2$ and define
$H_a'(u_1,v_1)=(p_1^*H_a)_{P_2}(u_1,v_1)/u_1^{m_1(H_a)}$,
which is a polynomial by Assertion (4.3.8).
By a similar discussion as above, we conclude:
$d(P_2;p_1^*f(u_1,v_1))=d(P_2;p_1^*H_a(u_1,v_1)^{n/a})
=d(P_2;p_1^* R_a(u_1,v_1))$,
$(p_1^*f)_{P_2}=(p_1^*H_a)_{P_2}^{n/a}+(p_1^* R_a)_{P_2}$ and
$\deg_{v_1}((v_1^{a_2}+\xi_2 u_1^{b_2})^{A_3}
-{H_a'}^{n/a}(u_1,v_1))\le A_2-a/a_1-1$.
In other words, $H_a'(u_1,v_1)$ is the
$n/a$-th \TS~ polynomial of $(v_1^{a_2}+\xi_2 u_1^{b_2})^{A_3}$.
In particular,  if $a_1a_2$ divides $a$, $A_3/(n/a)=a/a_1a_2$ is an integer and
we can see easily that
$H_a'(u_1,v_1)=(v_1^{a_2}+\xi_2 u_1^{b_2})^{a/a_1a_2}$ if $a_1a_2|a $ and
$p_1^*H_a(u_1,v_1)=u_1^{m_1(H_a)}(v_1^{a_2}+\xi_2 u_1^{b_2})^{a/a_1a_2}+
\text{(higher terms)}$.
Putting $h_2=H_{a_1a_2}$,  $C_2=D_{a_1a_2}$ and
 $a=a_1a_2$, we observe that
$p_1^*h_2(u_1,v_1)=u_1^{m_1(h_2)}(v_1^{a_2}+\xi_2 u_1^{b_2})+$
(higher terms) and therefore $p_1^*h_2(u_1,v_1)$
is clearly  non-degenerate.
\ninpar
(4.4) {\bf Inductive construction of a tower.}
Let
$\Cal T_j =\{
X_j\mapright{p_j} X_{j-1}\to\cdots \to X_1\mapright{p_1} X_0=\bC^2\}$
be a tower  of  toric modifications
with the corresponding weight vectors
$P_i={}^t(a_i,b_i)$
such that  $a_1\cdots a_j|n$  and $a_i\ge 2, ~i=1,\dots, j.$
Put
$A_{i+1}:=n/{a_1\cdots a_i}, i\le j$ and for simplicity
$h_i(x,y)=H_{a_1\cdots a_i}(x,y)$,
$C_i=D_{a_1\cdots a_i}$ and $\Phi_i=p_1\circ\cdots\circ p_i:X_i\to X_0.$
Let $D_a^{(\ell)}$ and $C_i^{(\ell)},~(i\ge \ell)$
be the strict transforms of $D_a$ and $C_i$ to $X_\ell$ respectively.
The map $p_{i}:X_{i}\to X_{i-1}$ is an admissible
 toric modification for $\Phi_{i-1}^*f$
associated with a regular simplicial cone subdivision
$\Si_{i-1}^*$.  Let
$\Xi_{i}=C^{(i)}\cap X_i$ be the
center of the  modification  $p_{i+1}$
and let $(u_i,v_i)$ be the chosen modification
local coordinate
system with the center $\Xi_i$  so that $\{u_i=0\}$ is the defining equation of
the exceptional divisor  $E_i:=\what E(P_i)$ for $i=1,\dots, j$. We assume  the
following properties (1-j), (2-j) and (3-j) for the tower.
\nl
{\bf (1-j)}  $(C_i,O)$ is a germ of an irreducible curve
at the origin for  $i=1,\dots, j$
and the strict transform $C_i^{(i)}$ to $X_i$ is smooth and is defined by
$\{v_i=0\}$. The pull backs of $f$ and  $h_i,~i\le j$
equal :
$$
\Phi_i^*f(u_i,v_i)=\cases
&u_j^{m_j(f)}v_j,~~\qquad i=j~ \text{and}~A_{j+1}=1\\
&u_i^{m_i(f)} (v_i^{a_{i+1}}+\xi_{i+1}u_i^{b_{i+1}})^{A_{i+2}}+
\text{(higher terms)},~~\quad \text{otherwise}\endcases
\tag 4.4.1$$
$$
\Phi_i^*h_\ell(u_i,v_i)=
\cases
&u_i^{m_i(h_i)}v_i,\qquad i=\ell\\
&u_i^{m_i(h_{\ell})}
(v_i^{a_{i+1}}+\xi_{i+1}u_i^{b_{i+1}})^{A_{i+2}/A_{\ell+1}}+
\text{(higher terms)},\qquad i<\ell
\endcases
 \tag 4.4.2
$$
The modification coordinates $(u_i,v_i)$ are characterized by (4.4.2).
We assume  $a_{j+1}\ge 2$ in (4.4.1), if $A_{j+1}\ge 2$.
More generally, for any positive integer
$a$ with $a|n$ and $a_1\cdots a_{i+1}|a$,
we have
$$
\Phi_i^*H_{a}(u_i,v_i)=u_i^{m_i(H_a)}
(v_i^{a_{i+1}}+\xi_{i+1}u_i^{b_{i+1}})^{a/a_1\cdots a_{i+1}}
+\text{(higher terms)}\tag 4.4.3$$
Here $m_i(h_i),~m_i(H_a)$ and $m_i(f)$ are the respective
multiplicities
of the pull backs $\Phi_i^*h_i,~\Phi_i^*H_a$  and  $\Phi_i^*f$
on the exceptional divisor $E_i$ and they satisfy the
equalities:
$$
m_i(h_i)\times A_{i+1}=m_i(H_a)\times n/a=m_i(f),~
m_1(f)= a_1b_1A_2,~
m_i(f)=a_i m_{i-1}(f)+a_ib_iA_{i+1}
$$
{\bf (2-j)} The local intersection multiplicities at the origin are given by
$$\cases
I(C_i,C;O)&=\sum_{s=1}^{i+1}a_sb_sA_{s+1}^2\big/A_{i+1},\\
I(C_i,C_\ell;O)&=\sum_{s=1}^{i+1}a_sb_sA_{s+1}^2 \big/(A_{i+1}A_{\ell+1}),
\quad 1\le i<\ell\le j
\endcases
$$
More generally, $I(D_a,C;O)=\sum_{s=1}^{j}a_sb_sA_{s+1}^2
\big/(n/a)+I(D_a^{(j)},C^{(j)};\Xi_j),$ if $a|n$ and $a_1\cdots a_j|a.$ \nl
{\bf (3-j)} For any non-zero polynomial $\al(x,y)\in \bC\{x\}[y]$
with $\deg_y \al(x,y)<a_1\cdots a_j$,
the pull back $\Phi_j^* \al$ can be written as
$\Phi_j^* \al=U\times u_j^s$  in a small neighbourhood
$W_j$ of $\Xi_j$ for some integer $s\ge 0$.
\par
If $A_{j+1}=1$, then $h_{j}=f$ and (4.4.1) says that $\Phi_j:X_j\to X_0$
is a good resolution of $C$.
If  $A_{j+1}\ge 2$,
we will add to the tower a toric
modification $p_{j+1}: X_{j+1}\to X_j$ keeping the above properties.
 Let $P_{j+1}={}^t(a_{j+1},b_{j+1})$ be
 the weight vector of the unique face of
$\Ga(\Phi_j^*f;(u_j,v_j))$ characterized by (4.4.1) and (4.4.2):
$
\Phi_j^*f(u_j,v_j)=
 u_j^{m_j(f)}
(v_j^{a_{j+1}}+\xi_{j+1}u_j^{b_{j+1}})^{A_{j+2}}+$
(higher terms).
Choose a  regular simplicial cone
subdivision $\Si_{j}^*$ of the
$\Ga^*(\Phi_{j}^*f;(u_j,v_j))$ and make  the corresponding
modification $p_{j+1}:X_{j+1}\to X_j$ with center $\Xi_j\in E_j$.
Then  $\Phi_j^*h_{j+1}(u_j,v_j)$ is non-degenerate by (4.4.2), so in the right
toric chart $\si=(P_{j+1},P_{j+1}')$
we can write
$\Phi_{j+1}^* h_{j+1}(x_{\si},y_\si)=x_\si^{m_{j+1}(h_{j+1})}
y_\si^{m_{j+1}'(h_{j+1})}\((y_\si+\xi_{j+1})+x_\si  G\)
$
where  $m_{j+1}(h_{j+1})$ and $m_{j+1}'(h_{j+1})$ are the
 multiplicities on
$E_{j+1}=\what E(P_{j+1})$ and $\what E(P_{j+1}')$ respectively.
The functions
$u_{j+1}:=x_\si$ and
 $v_{j+1}:=y_\si^{m_{j+1}'(h_{j+1})}\((y_\si+\xi_{j+1})+x_\si  G\)
$
give a  system of coordinates
 in a neighborhood $W_{j+1}$ of  the intersection point $\Xi_{j+1}$
of $C_{j+1}^{(j+1)}$ and $E_{j+1}$.
By the definition   the strict transform $C_{j+1}^{(j+1)}$
is smooth and  is defined by $\{v_{j+1}=0\}$ in $W_{j+1}$.
We show (3-(j+1)) first.
For $\al(x,y)\in \bC\{x\}[y]$ with
$\deg_y a< a_1\cdots a_{j+1},$
its  \Ts~expansion \wrt~ $h_j$:
$\al(x,y)=\sum_{i=1}^{a_{j+1}}\al_i(x,y) h_{j}^{a_{j+1}-i}(x,y)$ has
coefficients with $ \deg_{y}\al_i<a_1\cdots a_{j}$.
Applying (3-j) inductively if $\al_i\ne 0$ yields
$\Phi_j^*(\al_i h_{j}^{a_{j+1}-i})=U_i u_j^{\nu_i}v_j^{a_{j+1}-i}$
with $\nu_i \ge 0$ and a unit $U_i$.
So by Lemma (4.1),  $\Phi_{j+1}^* \al=p_{j+1}^*(\Phi_j^*\al)=U\times u_{j+1}^s$
for a unit $U$ on $W_{j+1}$ and  $s \ge 0$. \nl
\indent
If $a_{j+1}=A_{j+1}$ i.e., $A_{j+2}=1,$ the modification $\Phi_{j+1}:X_{j+1}\to
X_0$ is a good resolution of $C,$ so clearly we have
(1-(j+1)) and (2-(j+1)).
If  $A_{j+2}\ge 2,$ we write
$$
\Phi_{j+1}^*f(u_{j+1},v_{j+1})=
 u_{j+1}^{m_{j+1}(f)}
(v_{j+1}^{a_{j+2}}+\xi_{j+2}u_{j+1}^{b_{j+2}})^{A_{j+3}}+
\text{(higher terms)}\tag 4.4.4$$
Note that
$m_{j+1}(f)=a_{j+1}m_j(f)+a_{j+1}b_{j+1} A_{j+2}$.
Using the $A_{j+2}$-th \Ts~expansions of $f$:
$f(x,y)=h_{j+1}^{A_{j+2}}+\sum_{i=2}^{A_{j+2}}c_{j+1,i}h_{j+1}^{A_{j+2}-i},$
and repeating the argument in (4.3), will prove $a_{j+2}\ge 2.$
As above, if $c_{j+1,i}\ne 0,$ write
$\Phi_{j+1}^*(c_{j+1,i}h_{j+1}^{A_{j+1}-i})(u_{j+1},v_{j+1})=U_{{j+1},i}
u_{j+1}^{m_{i}}v_{j+1}^{A_{j+1}-i}$
for some integer $m_i$ and a unit $U_{{j+1},i}$.
The Newton principal part $\Cal N(\Phi_{j+1}^*f)(u_1,v_1)$
contains the exponent
$(m_{j+1}(f)+b_{j+2},A_{j+2}-a_{j+2})$ by (4.4.4) and
 we conclude
 that $a_{j+2}\ge 2$ as in (4.3).
\par
Now we show  (1-(j+1)).
For $a$ with $a|n$ and $a_1\cdots a_{j+1}|a,$ consider the \Ts~expansions:
$$
f(x,y)=
H_a^{n/a}+\sum_{i=2}^{n/a}c_{i}H_a^{n/a-i},\qquad
H_a=h_{j+1}^{\beta_{j+1}}+\sum_{i=2}^{\beta_{j+1}}
d_{i}h_{j+1}^{\beta_{j+1}-i}$$
with $\deg_y c_i<a$ and
$\deg_y d_i<a_1\cdots a_{j+1}$
where $\beta_i:=a/a_1\cdots a_i$.
Applying the same argument to the $h_{j+1}$-expansion of
$R:=f-H_a^{n/a}=\sum_{i=2}^{n/a}c_{i}H_a^{n/a-i}$,
we see:
$\deg_{v_{j+1}}\Phi_{j+1}^*(R)<A_{j+2}-\beta_{j+1}$.
But from (4.4.3) with $g_a,G_a \in \bC\{u_{j+1},v_{j+1}\}$ follows:
$$\align
\Phi_{j+1}^*H_a(u_{j+1},v_{j+1})&=p_{j+1}^*(\Phi_j^*H_a)(u_{j+1},v_{j+1})
=u_{j+1}^{m_{j+1}(H_a)}(v_{j+1}^{\beta_{j+1}}
+u_{j+1}g_a(u_{j+1},v_{j+1}))\\
\Phi_{j+1}^*H_a^{n/a}(u_{j+1},v_{j+1})&=
p_{j+1}^*(\Phi_j^*H_a^{n/a})(u_{j+1},v_{j+1})
=u_{j+1}^{m_{j+1}(f)}(v_{j+1}^{A_2}+u_{j+1}G_a(u_{j+1},v_{j+1}))
\endalign
$$
So
$B_a:=(m_{j+1}(H_a),\beta_{j+1})$ is the left  end vertex of $\Ga(\Phi_{j+1}^*
H_a),$
$n/a\times B_a$  is the left end vertex of $\Ga(\Phi_{j+1}^* H_a^{n/a})$
and also  of $\Ga(\Phi_{j+1}^*f;(u_{j+1},v_{j+1}))$.
By the arguments of (4.3) and (4.4.3), the first face $\De_a$ of
$\Ga(\Phi_{j+1}^* H_a;(u_{j+1},v_{j+1})),$
which contains $B_a,$ has the weight vector $P_{j+2}={}^t(a_{j+2},b_{j+2}),$
hence
$$\cases
d(P_{j+2};\Phi_{j+1}^*f)=d(P_{j+2};\Phi_{j+1}^*H_a^{n/a})=
d(P_{j+2};\Phi_{j+1}^*H_a)\times n/a&\\
\deg_{v_{j+1}}((\Phi_{j+1}^*f)_{P_{j+2}}-
(\Phi_{j+1}^*H_a)_{P_{j+2}}^{n/a})(u_{j+1},v_{j+1})<
A_{j+2}-\beta_{j+1}&
\endcases
\tag4.4.5
$$
Note:
 $\deg_{v_{j+1}}\Phi_{j+1}^*f=d(P_{j+1};\Phi_{j}^*f)$.
The polynomial $H_a'(u_{j+1},v_{j+1}):=(\Phi_{j+1}^* H_a)_{P_{j+2}} /
 u_{j+1}^{m_{j+1}(H_a)}$ is  monic in $v_{j+1}$  of degree
$\beta_{j+1}=a/a_1\cdots a_{j+1},$ implying
with the inequality of (4.4.5) the
\proclaim{Assertion (4.4.6)} If $a_1\cdots a_{j+1}|a,$ then
 $H_a'(u_{j+1},v_{j+1})$ is the
$n/a$-th \TS~
 polynomial of
$(\Phi_{j+1}^*f)_{P_{j+2}}(u_{j+1},v_{j+1})/u^{m_{j+1}(f)}
=(v_{j+1}^{a_{j+2}}+\xi_{i+1}u_{j+1}^{b_{j+2}})^{A_{j+2}} \in
\bC\{u_{j+1}\}[v_{j+1}].$ In particular, if
$a_1\cdots a_{j+2}|a$, then
$H_a'(u_{j+1},v_{j+1}) =
(v_{j+1}^{a_{j+2}}+\xi_{i+1}u_{j+1}^{b_{j+2}})^{\beta_{j+2}},
\Phi_{j+1}^*H_a(u_{j+1},v_{j+1})= u_{j+1}^{m_{j+1}(H_a)}
(v_{j+1}^{a_{j+2}}+\xi_{i+2}u_{j+1}^{b_{j+2}})^{\beta_{j+2}}+$
(higher terms), with $\beta_{j+2}:=a/a_1\cdots a_{j+2}$.
\endproclaim
This proves  (1-(j+1)).
The assertion about the intersection multiplicities
(2-(j+1)) follows immediately from  Lemma (3.4.2).
\nlind
As $a_1\cdots a_i$ divides $n$ and $a_i\ge 2$ for each
$i=1,\dots, k$,
the above inductive construction stops after a finite number of toric
modifications.  In fact,
 $k$ (respectively $k-1$)
is the number of Puiseux pairs if $b_1>1$ (resp. if $b_1=1$.)
 See  [27]  and [17].
Thus we have proved the following.
\proclaim{Theorem (4.5)}Let
 $f(x,y)\in\bC\{x\}[y]$ be monic of degree $n$ with the initial expansion
$$f(x,y)=(y^{a_1}+\xi_1 x^{b_1})^{A_2}+\text{(higher terms)},\quad
n=a_1A_2,\quad
a_1>1$$
and defining  in a neighbourhood $W_0$ an irreducible curve $C:=\{(x,y)\in
W_0;f(x,y)=0\}$   at the origin.
There exits a
resolution tower $\Cal T,$ satisfying the following conditions (1) and (2), of
toric modifications:
$\Cal T=\{~X_k\mapright{p_k} X_{k-1}\to\cdots
\to X_1\mapright{p_1} X_0=\bC^2\}$
having the  weight vectors
$\{P_i={}^t(a_i,b_i);i=1,\dots, k\}$
where $n=a_1\cdots a_k,~a_i\ge 2,~i=1,\dots,k.$
With $A_i=a_ia_{i+1}\cdots a_{k},$   let $h_i(x,y)$ be the
$A_{i+1}$-th \TS~ polynomial of $f(x,y)$
and let $C_i=\{(x,y)\in \bC^2;h_i(x,y)=0\},$ $i=1,\dots,k$.
Note $h_k=f$ and  $C_k=C$.
Denote by $\Xi_{i}\in E_i:=\what E(P_i)$   the center of
 $p_{i+1},$
by  $(u_i,v_i)$  the modification
 local coordinate
 centered at
 $\Xi_i$ so that $\{u_i=0\}$ is the defining equation of
the  divisor $E_i.$
Put
$\Phi_i=p_1\circ\cdots\circ p_i:X_i\to X_0$.
\nl
(1) For each $i=1,\dots, k$, $C_i$ is an irreducible curve at the origin
having  the good resolution $\Phi_i,$ such that
the strict transform $C_i^{(i)}$ in $X_i$ is  defined by
$\{v_i=0\}$.
The pull backs are
$$
\Phi_i^*h_{\ell}(u_i,v_i)=\cases
u_i^{m_i(h_i)}v_i,\qquad &i=\ell\\
u_i^{m_i(h_{\ell})}
(v_i^{a_{i+1}}+\xi_{i+1}u_i^{b_{i+1}})^{A_{i+2}/A_{\ell+1}}
+\text{(higher terms)},~&i<\ell\endcases
$$ In particular, putting $\ell=k$,
$$
\Phi_i^*f(u_i,v_i)=\cases
u_k^{m_k(f)}v_k,\quad &i=k\\
u_i^{m_i(f)}(v_i^{a_{i+1}}+\xi_{i+1}u_i^{b_{i+1}})^{A_{i+2}}+
\text{(higher terms)},
{}~&i<k
\endcases
\tag 4.5.1
$$
where  the
multiplicities $m_i(h_\ell)$ and $m_i(f)$
multiplicities
of the pull backs $\Phi_i^*h_\ell$  and  $\Phi_i^*f$
on $E_i$  satisfy the
equalities:
$
m_i(h_\ell)=m_i(f)/A_{\ell+1}$, $m_1(f)= a_1b_1A_2$ and
$m_i(f)=a_i m_{i-1}(f)+a_ib_iA_{i+1}$
for $i=1,\dots, \ell$.
More explicitly
$$\cases
m_i(f)&=a_ib_iA_{i+1}+a_ia_{i-1}b_{i-1}A_i+\cdots+a_i\cdots a_1b_1A_2
=(\sum_{\ell=1}^i a_\ell b_\ell A_{\ell+1}^2)/A_{i+1}\\
m_i(h_\ell)&=(\sum_{j =1}^i a_j  b_j  A_{j +1}^2)/
(A_{i+1}A_{\ell+1}),\quad \qquad\qquad i\le \ell
\endcases\tag 4.5.2
$$
(2) The local intersection multiplicities
are
$I(C_\ell,C_{j};O)=\sum_{i=1}^{\ell+1}a_ib_i A_{i+1}^2/(A_{\ell+1}A_{j+1}),
\ell<j\le k$.
\endproclaim
The equality (4.5.2) follows from  (4.5.1).
The other assertions are etablished in the inductive argument.
\definition{Definition (4.5.3)}
The toric tower  of Theorem (4.5) is
{\it a \Ts~  resolution tower of toric modifications} of
 $C$, the coordinates $(u_i,v_i)$ of $W_i$ are
{\it  \Ts~ coordinates centered
at} $\Xi_i$, and the curve $C_i$ is
{\it the $A_{i+1}$-th \TS~ curve} of $C$.
\enddefinition
The combinatorial choice of the admissible subdivisions $\Si_i^*$'s
determins completely the \nl \Ts~  resolution tower of toric modification.
In Theorem (4.7), we will show that the length of the
tower $k$ and the sequence of
the weight vectors $\{P_1,\dots,P_k\}$ are independent of the choice of
a certain resolution tower of toric modifications.
\ninpar
{\bf B. Intersections of other \TS~ polynomials.}
Let, as before,  $H_a(x,y)$ be the $n/a$-th
\TS~ polynomials and
$D_a=\{H_a(x,y)=0\}$.
\comment
The following theorem describes the relation of
the \Ts~ curves
$\{C_i; i=1,\dots k\}$ which we have already studied
and
other \Ts~
curves $\{D_a;a\ne a_1\cdots a_i\}$.\endcomment
\proclaim{Theorem (4.6)}
If $a|n, ~a_1\cdots a_s|a$, $ a_1\cdots a_{s+1}\not| a$
and $a\ne a_1\cdots a_s,$
then $D_a$ and $C$
have the same toric tangential direction of depth  $s$  and
$I(D_a,C_i;O)= \sum_{j=1}^{\al+1} a_j b_j A_{j+1}^2/(A_{i+1} n/a)$
where $ \al=\min(s,i)$.
\endproclaim
 \demo{Proof}
Recall that
$\Phi_s^*f(u_s,v_s)=u_s^{m_s(f)}(v_s^{a_{s+1}}
+\xi_{s+1}u_s^{b_{s+1}})^{A_{s+2}}+\text{(higher terms)}$.
We consider the face function of the pull-back $\Psi_s^*H_a$ and put
$H_a'(u_s,v_s):=(\Psi_s^*H_a)_{P_{s+1}}(u_{s},v_s)/u_s^{m_s(H_a)}.$
We have seen in the inductive construction of the
\Ts~ tower that $D_a$ has the same toric tangential direction at least
 of depth
 $s$ with $C$.
We have shown in Assertion (4.4.18) that $H_a'(u_s,v_s)$ is the $n/a$-th
\TS~ polynomial of $(v_s^{a_{s+1}}+\xi_{s+1}u_s^{b_{s+1}})^{A_{s+2}}$.
Now the main step of the proof is the  following.
\proclaim{Lemma (4.6.1)}
 The constant term of the polynomial $H_a'(u_s,v_s) \in \bC\{u_s\}[v_s]$ is
zero  and
$v_s^{a_{s+1}}+\xi_{s+1}u_s^{b_{s+1}}$ does not divide
$H_a'(u_s,v_s)$.
\endproclaim
\demo{Proof}
Put
$ \beta_j=a/a_1\cdots a_j$.
 The point is that $\beta_{s+1}:=A_{s+2}/(n/a)$ is not an integer.
 As  $H_a'(u_s,v_s)$ is the $n/a$-th
\TS~ polynomial of $(v_s^{a_{s+1}}+\xi_{s+1}u_s^{b_{s+1}})^{A_{s+2}}$,
we have
$H_a'(u_s,v_s)=v_s^{\beta_s}\Trunc^{([\beta_{s+1}])}
(1+\xi_{s+1}u_s^{b_{s+1}}v_s^{-a_{s+1}})^{\beta_{s+1}}
=v_s^{\beta_s}\sum_{j=0}^{[\beta_{s+1}]}
 \binom{\beta_{s+1}}{j}
(\xi_{s+1}u_s^{b_{s+1}}v_s^{-a_{s+1}})^j$ by Lemma (2.3).
Thus $H_a'(u_s,v_s)$ does not have a constant term as a polynomial
of $v_s$.
If
$v_s^{a_{s+1}}+\xi_{s+1}u_s^{b_{s+1}}$  divides
$H_a'(u_s,v_s),$ we will get a contradiction. In fact, the polynomial
$$H_a'(u_s,v_s)=(v_s^{a_{s+1}}+\xi_{s+1}u_s^{b_{s+1}})H_a''(u_s,v_s)
\tag 4.6.4$$
is the $n/a$-th \TS~polynomial of
$(v_s^{a_{s+1}}+\xi_{s+1}u_s^{b_{s+1}})^{A_{s+2}-n/a}$.
By the generalized binomial formula again, we have
$H_a''(u_s,v_s)=
v_s^{\beta_s-a_{s+1}}\sum_{j=0}^{[\beta_{s+1}]-1}
 \binom{\beta_{s+1}-1}{j}
(\xi_{s+1}u_s^{b_{s+1}}v_s^{-a_{s+1}})^j$.
Comparing the coefficients of
$v_{s+1}^{\beta_s-[\beta_{s+1}]a_{s+1}}$
in (4.6.4),
we get:
$\binom{ \beta_{s+1}}{[\beta_{s+1}]}=
\binom{\beta_{s+1}-1}{[\beta_{s+1}]-1},$ which is a contradiction
as $\beta_{s+1}\ne [\beta_{s+1}].$  ~\QED
\enddemo
Now by  the Lemma the curve $D_a$ has  the same toric tangential direction
 of depth  $s$ but not of depth $s+1$ with $C$. In particular, $D_a^{(s+1)}\cap
C^{(s+1)}=\emptyset$.
 The main problem in  proving the assertion about the intersection
multiplicity is that $D_a$ may be neither  irreducible nor reduced.
See Example (4.9).
Let $D_{a,1},\cdots, D_{a,\ell}$ be the irreducible components.
Let $k_a(u_s,v_s)$ and $k_{a,j}(u_s,v_s),j=1,\dots,\ell$
be the defining functions of the strict transforms
$D_a^{(s)}$ and $D_{a,j}^{(s)}$ for $j=1,\dots,\ell$.
Then we can write
$
k_a(u_s,v_s)=v_s^{\beta_s}+\sum_{t=1}^{\beta_s}\ga_t(u_s)v_s^{\beta_s-t}$,
$\ga_t(u_s)\in \bC\{u_s\}$ and
$$k_{a,j}(u_s,v_s)=\cases
(v_s^{a_{s+1,j}}+\xi_{a,j}u_s^{b_{s+1,j}})^{A_{s+2,j}}
+\text{(higher terms)},
{}~ & b_{s+1,j}\ne 0\\
v_s^{A_{s+2,j}}U_j, \phantom{aaaaaaaaaaaaaaaaa}b_{s+1,j}=0,~&
a_{s+1,j}=1\endcases\tag 4.6.5
$$
where $\gcd(a_{s+1,j},b_{s+1,j})=1$ and $U_j$ is a unit. They satisfy:
$$\beta_s=\sum_{j=1}^\ell a_{s+1,j}A_{s+2,j}\tag 4.6.6$$
Recall that the weight vector of the unique face of $\Ga(k_{a,j};(u_s,v_s))$
corresponds to the weight vector of a face of $\Ga(k_a;(u_s,v_s))$.
By Assertion (4.4.18),
 the Newton boundary $\Ga(k_a;(u_s,v_s))$ starts  with the
face (possibly a vertex) of the weight vector $P_{s+1}$ and
any other face has a milder slope.
Therefore we have
$
b_{s+1,j}/a_{s+1,j}\ge b_{s+1}/a_{s+1}$ if
$b_{s+1,j}\ne 0$.
Now we apply Lemma (3.4.2) to compute the intersection numbers.
For $i\le s$, we have
$
I(D_a,C_i)=\sum_{j=1}^{i+1} a_jb_j A_{j+1}^2\big /(A_{i+1} n/a)$,
$i\le s$
and for $i\ge s$, with $P_{s+1,t}:={}^t(a_{s+1,t},b_{s+1,t})$
we have
$$\align
I(D_a,C_i)&=\sum_{j=1}^{s} a_jb_j A_{j+1}^2\big/
(A_{i+1}n/a)+\sum_{t=1}^\ell I(P_{s+1},P_{t,s+1})
A_{s+2,t}{A_{s+2}}\big/{A_{i+1}}\\
&= \sum_{j=1}^{s} a_jb_j {A_{j+1}^2}\big/({A_{i+1}}n/a)
+\sum_{t=1}^\ell b_{s+1}\beta_{s,t}
A_{s+2,t}{A_{s+2}}\big /{A_{i+1}}\quad \text{by (4.6.6)}\\
&=\sum_{j=1}^{s+1} a_jb_j {A_{j+1}^2}\big/({A_{i+1}}n/a),\quad i>s
\quad \text{by (4.6.5).}\qquad\qquad{\roman \QED}
\endalign
$$
\enddemo
\ninpar
{\bf C. Relations with other toric towers}.
 Consider two toric resolution towers:
$$
\Cal T=\{
X_k\mapright{p_k} X_{k-1}\to\cdots \to X_1\mapright{p_1} X_0=\bC^2\},\quad
\Cal Q=\{Y_{s}\mapright{q_{s}}Y_{{s}-1}\to \cdots \to
Y_1\mapright{q_1}Y_0=\bC^2\}$$
where $\Cal T$ is a Tschirnhausen tower of resolution
with the  weight vectors
$P_i={}^t(a_i,b_i)$, $i=1,\dots, k$ and
 $n=a_1\cdots a_k,~a_i\ge 2,~i=1,\dots,k$ as in Theorem (4.5).
Let $A_i=a_i a_{i+1}\cdots a_{s}$ and  let
 $h_{i}(x,y)$ be the
$A_{i+1}$-th \TS~ polynomial of $f(x,y)$
and let $C_{i}$ be the corresponding  \Ts~ curve
for $i=1,\dots,{k}$ as before. Let
$Q_i={}^t(\al_i,\beta_i)$, $i=1,\dots, {s}$ be
the corresponding weight vectors of $\Cal Q$ with
 $n=\al_1\cdots \al_{s}$.
We assume that $\al_i\ge 2,~i=1,\dots,{s}$
and  $Q_1=P_1$. We call such a toric tower $\Cal Q$
{\it a \Ts-good resolution tower}.
\comment
(In the case of $b_1=1$, $(x,y)$ are not  good coordinates in the sense of
[17] but we do not take any coordinate change
as we want  $f(x,y)$ to be
  monic.)
\endcomment
 A \Ts~ resolution tower is a \Ts-good resolution tower
by Theorem (4.5).
Now Theorem (4.5) can be generalized as follows.
\proclaim{Theorem (4.7)}Let $f(x,y)$ as in Theorem (4.5).
Let $\Cal T$ and $\Cal Q$ be  as above.
Assume that $q_{i+1}:Y_{i+1}\to Y_{i}$ is a toric modification centered at
$\Theta_{i}\in E_i':=\what E(Q_i)$ with
 the modification
 local coordinate system $(w_i,z_i),$
 so that $\{w_i=0\}$  defines
the  divisor $E_i'$.
Put
$\Psi_i=q_1\circ\cdots\circ q_i:Y_i\to Y_0$.
Then we have the following properties.\nl
{\rm (1)}(Uniqueness of the weight vectors) ${s}=k$ and $Q_i=P_i$ for
$i=1,\dots,k$.\nl
{\rm (2)} For each $i=1,\dots, {s}$,
 $\Psi_i :~Y_i\to Y_0$ gives a good resolution
of $C_{i}$ and
the pull backs of the polynomials are written (up to a multiplication
of a non-zero constant) as
$$
\Psi_i^*h_{\ell}(w_i,z_i)=\cases
{w_i}^{m_i'(h_{\ell})}
({z_i}^{a_{i+1}}+\theta_{i+1}{w_i}^{b_{i+1}})^
{A_{i+2}/A_{\ell+1}}
+\text{(higher terms)},~&i<\ell\\
{w_i}^{m_i'(h_{i})}{z_i'}
,\qquad\qquad\qquad &i={\ell}\endcases\tag 4.7.1
$$
where
$z_i'$ is either $z_iU_i$ with a unit $U_i$  or
$c_i((z_i+\eta_{i} w_i^{\ga_i})+\text{(higher terms)})$ with
$c_i,\eta_i\in \bC^*$
for some integer $\ga_i,~\ga_i>b_{i+1}/a_{i+1}$.
In particular, putting $\ell={s}$, we have
$$
\Psi_i^*f({w_i},{z_i})=\cases
({z_i}^{a_{i+1}}+\theta_{i+1}{w_i}^{b_{i+1}})^{A_{i+2}}
+\text{(higher terms)},
\quad &i<{s}\\
{w_{s}}^{m_{s}'(f)}z_{s}',\qquad\qquad\qquad &i=s
\endcases\tag4.7.2
$$
where the
multiplicities
$m_i'(h_\ell')$ and $m_i'(f)$
of the pull backs ${\Psi_i}^*h_{\ell}$  and  ${\Psi_i}^*f$
on $E_i'$  satisfy the same inductive
equalities:
$$\cases
m_i'(h_{s})&=m_i'(f)/A_{{s}+1},\quad i\le {s}\\
m_1'(f)&=a_1b_1A_2~\quad
m_i'(f)=a_i m_{i-1}'(f)+a_ib_iA_{i+1}
\endcases\tag 4.7.1
$$
Thus we have also the uniqueness of the multiplities:
 $m_i'(h_s)=m_i(h_s)$ and $m_i'(f)=m_i(f)$.
\endproclaim
\demo{ Proof}
We consider the tower $\Cal Q$.
Let $\wtl \al_1=\min(\al_1,\beta_1) $ and $\wtl \beta_1=\max(\al_1,\beta_1)$
and let
 $n_1=\wtl \al_1,~m_1=\wtl \beta_1$
and  $ n_i= \al_i$,~
$ m_i=\beta_i+\beta_{i-1}\al_i+\cdots+  \beta_2\al_3\cdots \al_i+ \wtl \beta_1
\al_2\cdots \al_i$
for  $i\ge 2$.
Then we have shown in Corollary (6.8) of [27] that
the Puiseux pairs of $C$ is
 given by
$\{(n_i,m_i);~i=1,\dots,s\},~(\beta_1>1)$ or
$\{(n_i,m_i);~i=2,\dots,s\},~(\beta_1=1)$.
 The same assertion is true
for the Tschirnhausen tower $\Cal T$.
By the assumption $Q_1=P_1$ and by the uniqueness of the  Puiseux pairs,
we conclude that
$s=k$ and $Q_i=P_i$.
The assertion (4.7.1) for $\ell>i$ follows easily by the induction on $i$.
In fact, we know that $C_\ell$ is irreducible and
$I(C_\ell,C;O)=\sum_{s=1}^{\ell+1}a_sb_sA_{s+1}^2/A_{\ell+1}$.
So by Lemma (3.4.2), $C_\ell$ can not be separated from $C$ on
$Y_i,~i<\ell$. Thus we have the expression (4.7.1).
As $\Psi_{i-1}^*h_{i}(w_i,z_i)$
is non-degenerate, we can write
$${\Psi_i}^*h_{i}(w_i,z_i)=\cases
&c_i((z_i+\eta_{i} w_i^{\ga_i})+\text{(higher terms)}
),\quad c_i,~\eta_i\in \bC^*\quad\text{or}\\
&z_iU_i,\quad U_i:\text{a unit}
\endcases$$
Assume the first case. The formula for the intersection multiplicity
in Theorem (4.5)
says that
$I(C_i^{(i)},C^{(i)};\Theta_i)=b_{i+1}A_{i+2}$.
This is the case if and only if
$\ga_i\ge b_{i+1}/a_{i+1}$. As $b_{i+1}/a_{i+1}$
is not an integer, we have $\ga_i > b_{i+1}/a_{i+1}$.~\QED
\enddemo
\remark{\bf Remark (4.8)}
Theorem (4.7) can be proved without using the uniqueness of the Puiseux pairs
by comparision stage by stage of the formulae for the intersections for the two
towers.
\endremark
\example{Example (4.9)}Put
$f(x,y)=(y^4+x^3)^6+x^{17} y^3$.
The first toric modification $p_1: X_1\to X_0$ can be defined by
the subdivision
$$
\Si_0^*=\{P_{0,0},\dots, P_{0,5}\}=\left\{
\pmatrix 1\\ 0\endpmatrix,
\pmatrix 2\\ 1\endpmatrix,\pmatrix 3\\ 2\endpmatrix,
\pmatrix 4\\ 3\endpmatrix,
\pmatrix 1\\ 1\endpmatrix,
\pmatrix 0\\ 1\endpmatrix\right\}$$
with weight vector $ P_1=P_{0,3}$.
Let $\si_3=\Cone (P_{0,3},P_{0,4})$.
On the chart $\bC_{\si_3}^2$, we take $u_1=x_{\si_3}$ and $v_1=y_{\si_3}+1$.
Then $C^{(1)}$ is defined by
$\{(u_1,v_1)\in W_1; v_1^6+u_1^5+\text{(higher terms)}=0\}$.
Thus we need one more toric modification
$p_2: X_2\to X_1$ and we choose the modification with respect to
 $$\Si_1^*=\{P_{1,0};,\dots, P_{1,7}\}=\{\pmatrix 1\\ 0\endpmatrix,
\pmatrix 2\\ 1\endpmatrix,
\pmatrix 3\\ 2\endpmatrix,
\pmatrix 4\\ 3\endpmatrix,\pmatrix 5\\ 4\endpmatrix,
\pmatrix 6\\ 5\endpmatrix,
\pmatrix 1\\ 1\endpmatrix,\pmatrix 0\\ 1\endpmatrix\}
$$
with weight vector $P_2=P_{1,5}$.
The weight vectors of the tower are
$P_1={}^t(4, 3)$ and $ P_2={}^t(6, 5)$.
By  computation, we have
$n=24$ and the various \TS~ polynomials are:
$H_2(x,y)= y^2$, $ H_3(x,y)=y^3$,
$H_4(x,y)=h_1(x,y)=y^4+x^3$,
$H_6(x,y)=y^6+3/2 x^3 y^2$,
$H_8(x,y)=(y^4+x^3)^2$ and
$H_{12}(x,y)= (y^4+x^3)^3$.
The intersection multiplicities are given by
$I(D_a,C;O)=36$, 54, 77, 108, 154, 231
respectively for $a=2,3,4,6,8,12$.
This example shows that $D_a$ which is different from $C_i,1\le i\le k$
 is not necessarily irreducible or reduced.
The zeta-function and the Milnor number are given by Theorem (5.1) in \S 5:
$\zeta(t)=(1-t^{72})(1-t^{462})/(1-t^{24})(1-t^{18})(1-t^{77})$ and
$\mu(f)=416$.
\endexample
\ninpar
\S 5. {\bf The zeta-function of the monodromy.}
 Let $f(x,y)$ be a monic polynomial
in $y$ of degree $n$ and irreducible at the origin. Let
$\Cal T=\{X_k\mapright{p_k}X_{k-1}\to \cdots\to X_1
\mapright{p_1} X_0\}$
be a \Ts-good toric resolution tower
with the weight
vectors  $\{P_i={}^t(a_i,b_i);i=1,\dots, k\}$.
 We will describe
the invariants of the monodromy [1],[19] of $f$
at the origin using the data of
the \Ts-good resolution tower.
\nlind
 Let
$\Si_i^*$ be the regular simplicial cone subdivision which is used to
 construct the modification
$p_{i+1}:X_{i+1}\to X_i$ and let $\{P_{i,0},P_{i,1},\dots, P_{i,r_i}, P_{i,r_i
+1}\}$
be the vertices of $\Si_i^*$ so that
$P_{i,0}={}^t(1,0)$ and $P_{i,r_i+1}={}^t(0,1)$.
Let $P_{i,j}={}^t(a_{i,j},b_{i,j})$. We assume that $P_{i+1}=P_{i,n_i}$
for $i=0,\dots, k-1$.
Note that, as $\det(P_{i,0},P_{i,1})=\det(P_{i,r_i},P_{i,r_i+1})=1$,
$P_{i,1}$ and
 $P_{i,r_i}$ have the forms $P_{i,1}={}^t(a_{i,1},1)$
and  $P_{i,r_i}={}^t(1,b_{i,r_i})$ respectively. This  implies that
  $n_i<r_i$.
The  configuration of the  exceptional divisors
$\{\what E(P_{i,j});j=1,\dots, r_i\}$
is a line configuration and $\what E(P_{i,0})$
is nothing but $\what E(P_{i-1,n_{i-1}})$.
Thus the exceptional divisors of the resolution
$\Phi_k:X_k\to X_0$ is the union of the strict transforms
$\{\what E(P_{i,j});0\le i\le k-1, 1\le j\le r_i\}$.
Let $m_{i,j}$ be the multiplicity of the pull-back
$\Phi_{i+1}^*f$ along $\what E(P_{i,j})$ and
let $\de_{i,j}$ be the number of irreducible components
of the divisor $(\Phi_k^*f)$ which intersect with $\what E(P_{i,j})$.
By Theorem 3 of  [1], the zeta-function $\zeta(t;O)$ of the  monodromy of
$f(x,y)$ is determined by those $\what E(P_{i,j})$ with $\de_{i,j}\ne 2$.
As we have seen in \S 3, $m_{i,j}=d(P_{i,j};\Phi_i^*f)$
and
 $$
\de_{i,j}=\cases 3 \quad &j=n_i\\
                1\quad &j=r_i\\
                2\quad &\text{otherwise}\endcases
\quad i\ge 1,\quad
\de_{0,j}=\cases  3 \quad &j=n_1\\
                1\quad &j=1~\text{or}~r_1\\
                2\quad &\text{otherwise}\endcases
,\quad i=0
$$
If  $n_0=1$, we subdivide $\Cone (P_{0,0},P_{0,1})$ so that
we can assume that $n_0>1$. Note that $\de_{k-1,n_{k-1}}=3$ as
$\what E(P_{k-1,n_{k-1}})=E(P_k)$ and it  intersects
 with $C^{(k)}$.
Recall that  the multiplicity $m_{i,n_i}$ is given by
$m_{i,n_i}=d(P_{i+1};\Phi_{i}^*f)=m_{i+1}(f)=a_{i+1}m_i(f)+a_{i+1}b_{i+1}A_{i+2}$
in the same notation as in \S 4.
Thus we need  determine $m_{0,1}, m_{i,r_i}$ for $i=1,\dots,k$.
To determine $m_{i,r_i}$, we consider the expression by (4.7.2):
$\Phi_i^*f (u_i,v_i)=u_i^{m_i(f)}
(v_i^{a_{i+1}}+\xi_{i+1}u_i^{b_{i+1}})^{A_{i+2}}+\text{(higher terms)}$
for $ i<k$.
As $\Si_i^*$ is assumed to be admissible for $\Phi_i^*f$, we know that
$(m_i(f)+b_{i+1}A_{i+2},0)\in\De(P_{i,r_i};\Phi_i^*f)$.
This  observation and the expression  $P_{i,n_i}={}^t(1,b_{i,n_i})$
implies  that
$m_{i,r_i }=m_i(f)+b_{i+1}A_{i+2}=m_{i+1}(f)/a_{i+1}$.
Finally as $(0,n)\in \De(P_{0,1};f)$,
 we have  that $m_{0,1}=A_1$  by a similar argument as above.
Thus applying  Theorem 2 of [1], we obtain the first part of
\proclaim{Theorem (5.1)} The zeta-function and Milnor number of  $f(x,y)$ are:
$$
\zeta(t;O)=\frac 1{(1-t^{A_1})}
\prod_{i=1}^k\frac{(1-t^{m_i(f)})}{(1-t^{m_i(f)/a_i})},
\quad
\mu(f;O) =1-A_1+\sum_{i=1}^k(A_{i}-1) b_iA_{i+1}
 $$
\endproclaim
\demo{Proof}
 By  the equality $-1+\mu(f;O)=\deg \zeta(t;O)$, we have
$$\align
-1+\mu(f;O)=
&=-A_1+\sum_{i=1}^k \(1-\dsize\frac1{a_i}\)m_i(f)\\
&=-A_1+\sum_{i=1}^k \(1-\dsize\frac1{a_i}\)
\(\sum_{\ell=1}^i a_\ell b_\ell A_{\ell+1}^2\)/A_{i+1}\\
&=-A_1+\sum_{\ell=1}^ka_\ell b_\ell A_{\ell+1}^2
\sum_{i=\ell}^k\(1-\dsize\frac1{a_i}\)/{A_{i+1}}\\
&=-A_1+\sum_{\ell =1}^k\(A_{\ell }-1\) b_\ell A_{\ell +1}.
\qquad \qquad \text{\QED}
\endalign
$$
\enddemo
\ninpar
\S 6. {\bf Conditions implying equi-singularity.}
Let $f_t(x,y)=f(x,y,t)\in \bC\{x,t\}[y]$ be an analytic family of monic
polynomial in $\bC\{x\}[y]$ of degree $n$  in $y$  defined for $t$ in an open
connected neighborhood
$U$  of the origin in $\bC$.
Let
$C(t):=\{f_t(x,y)=0\},\quad t\in U,$
be the corresponding family of germs of curves at the origin.
We assume that $C(0)$ is irreducible and reduced at the origin
and that $f_t(x,y)$ has an initial expansion
$$
f_t(x,y)=(y^{a_1}+\xi_1 x^{b_1})^{A_2}+\text{(higher terms)}
$$
with $\xi_1 \ne 0$  independent of t  and $a_1\ge 2$.
Let
$X_k\mapright{p_k}X_{k-1}\to \cdots\to X_1\mapright{p_1}X_0=\bC^2$
be the \TS~ resolution tower of $(C(0),O)$  with the weight vectors
$\{P_i={}^t(a_i,b_i);i=1,\dots, k\}$.
We assume further that the $A_{i+1}$-th \TS~ polynomials $h_i(x,y,t)$ of
$f_t(x,y)$ for $i=1,\dots, k-1$ are independent of the parameter $t$. Note that
this is the case if the coefficients of $y^j$ do not depend on $t$
for any $ j\ge n-a_1\cdots a_{k-1}$. Consider the germs of curves
$C_i:=\{ h_i(x,y):=h_i(x,y,t)=0\},~\quad i=1,\dots, k-1.$
Finally we  assume that the  local
intersection multiplicities satisfy the inequalities:
$$
a_k I(C_{k-1},C(t);O)\le I(C(0),C(s);O) < +\infty,\quad
 \text{for any}~t, s,~\text{with}~ s\ne 0.\tag 6.1
$$
\proclaim{Theorem (6.2)} Under the above assumptions for the family $f_t(x,y)$
,
 the germes $C(t),t\in U,$ are irreducible
at the origin and have the same toric tangential directions of  depth
$k', ~k' \ge k-1$ .
The family of germs of plane curves
$\{(C(t),O);t\in U\}$
is an equi-singular family and $\Phi_k:X_k\to X_0$ gives a
simultaneous resolution for the family $\{C(t);t\in U\}$
where  $\Phi_k=p_1\circ\cdots\circ p_k$.
In particular, the Milnor number $\mu(f_t;O)$ is constant and
coincides with $\mu(f_0;O)$. Moreover,
if equality holds in $
a_k I(C_{k-1},C(t);O)\le I(C(0),C(s);O) \quad
 \text{for any}~t, s,~\text{with}~ s\ne 0,$ the  germes
 $C(t),t\in U$
 do not have the
same toric tangential direction of depth $k$.
\endproclaim
\demo{Proof}We fix $\tau\ne 0, \tau\in U$.
We  first assume that
  $C(\tau)$ is irreducible. The irreducibility will be proved later.
 Assume that $C(\tau)$
has the same toric  tangential direction with $C(0)$ of depth $\theta$,
$\theta\le k$.  Then
we can  write  $C^{(j)}(\tau)$ as
$$\cases
C^{(j)}(\tau)&=\{(u_j,v_j)\in W_j;f_t^{(j)}(u_j,v_j)=0\}\\
f_t^{(j)}(u_j,v_j)&= (v_j^{a_{j+1}'}+\xi_{j+1}'u_j^{b_{j+1}'})^{A_{j+2}'}+
\text{(higher terms)}
\endcases\tag6.2.1
$$
where
$A_{i+1}':=a_{i+1}'\cdots  a_{\theta+1}'A_{\theta+2}'$ for
$ i\le\theta$.
Let  $P_{\theta+1}':={}^t (a_{\theta+1}',b_{\theta+1}')$.
$P_{\theta+1}'$
is a primitive weight vector and
if $P_{\theta+1}'=P_{\theta+1}$, we must have
 $\xi_{\theta+1}'\ne \xi_{\theta+1}$ by the assumption.
Comparing (6.1.1) and (6.2.1) and by the assumption, we have
$\xi_i'=\xi_i$, $a_{i}'=a_{i}$, $ b_{i}'=b_{i}$ and
$A_{i+1}'=A_{i+1}$.
Assume first that $\theta\le k-1$.
By Lemma (3.4.2) and (6.2.2),
the local intersection multiplicity is given by
$$\align
I(C(\tau),C(0);O)&=\sum_{i=1}^{\theta}
a_ib_iA_{i+1}^2+I(P_{\theta+1},P_{\theta+1}')
A_{\theta+2}A_{\theta+2}'\\
&\qquad\quad\le\sum_{i=1}^{\theta}b_iA_{i+1}A_i'+b_{\theta+1}A_{\theta+2}A_{\theta+1}'
=\sum_{i=1}^{\theta+1}b_ia_{i}A_{i+1}^2
\tag 6.2.3
\endalign
$$
where the equality holds if and only if
$a_{\theta+1}'b_{\theta+1}\le a_{\theta+1}b_{\theta+1}'$ or
$b_{\theta+1}'=0$.
On the other hand, by Theorem (4.5)  we have the equality:
$a_k I(C_{k-1},C(0);O)=\sum_{i=1}^{k}a_ib_iA_{i+1}^2.$
Thus (6.2.2) and (6.2.3) and the assumption (6.1.4)
 implies that we must have $\theta= k-1$ and
$ a _{k}'b_{k}\le a_{k}b_{k}'$ or $b_{k}'=0$ and
$I(C(\tau),C(0);O)=\sum_{i=1}^{k}a_ib_iA_{i+1}^2$.
We assert furthermore
$$b_{k}'\ne 0,\quad a_{k}'b_{k}= a_{k}b_{k}'.\tag6.2.4$$
In fact, assume first that $b_{k}'=0$.  Then
$C(\tau)=a_k C_{k-1}$ and $C(\tau)$ is not reduced.
 This is a contradiction to the assumption
$\dim_{\bC}\bC\{x,y\}/(f_\tau, h_{k-1})<\infty$.
Assume that $b_{k}'\ne 0$ and  $a_{k}'b_{k}<a_{k}b_{k}'$.
Then we get a contradiction:
$$\align
a_k I(C(\tau),C_{k-1};O)&=\sum_{i=1}^{k-1}
a_ib_iA_{i+1}^2+ a_k I(C_{k-1}^{(k-1)},C^{(k-1)}(\tau);O)
=\sum_{i=1}^{k-1}a_ib_iA_{i+1}^2+
a_kb_{k}'A_{k+1}'\\
&\qquad\qquad >\sum_{i=1}^{k-1}a_ib_iA_{i+1}^2+
a_k'b_kA_{k+1}'
=\sum_{i=1}^{k}a_ib_iA_{i+1}^2= I(C(\tau),C(0);O).
\endalign$$
Thus we have proved (6.2.4).
As $\gcd(a_k,b_k)=\gcd(a_k',b_k')=1$, (6.2.4) implies
$P_{k}'=P_k$ and  $A_{k+1}'=1$.
This also shows that $C^{(k)}(\tau)$ is smooth.
Thus under the assumption that $C(\tau)$ is irreducible at
the origin, we have proved  that $C(\tau)$ is reduced and
$\theta\ge k-1,\quad P_{k}'=P_k.$
This implies that $\mu(f_\tau;O)=\mu(f_0;O)$ by
applying  Theorem (4.5) to $C(\tau)$.
Note that  $\Phi_k:X_k\to X_0$
gives a simultaneous resolution of the family $\{C(\tau);\tau \in U\}$.
If $\theta=k$, the assertion is obvious and
$C^{(k)}(\tau)$
intersects with $C^{(k)}(0)$ at $\Xi_k$ and therefore
$I(C(0),C(\tau);O)>\sum_{i=1}^k a_ib_i A_{i+1}^2$.
This implies that the strict inequality in (6.1.3) must hold.
\par\noindent
{\bf Irreducibility of $\bold{C(\tau)}$}.
Now we prove that $C(\tau)$ is irreducible for any $\tau$.
Fix a $\tau$ and assume that $C(\tau)$ has $s$ irreducible components
at the origin and $s\ge 2$. Let $C(\tau;1),\dots, C(\tau;s)$ be the
irreducible components and let $C^{(j)}(\tau;1),\dots, C^{(j)}(\tau;s)$
be their strict transforms on $X_j$.
We assume that $C(\tau; i)$ has the same toric  tangential
direction of depth $\theta_i$ with $C(0)$. Then we can write
$$\cases
C^{(j)}(\tau;i)&=\{(u_j,v_j)\in W_j; f_{\tau,i}^{(j)}(u_j,v_j)=0\}\\
f_{\tau,i}^{(j)}(u_j,v_j)&=
 (v_j^{a_{i,j+1}}+\xi_{i,j+1}u_j^{b_{i,j+1}})^{A_{i,j+2}}+
\text{(higher terms)}, \quad j\le \theta_i\endcases
$$
where $A_{i,j}=a_{i,j}\cdots a_{i,\theta_i+1}A_{i,\theta_i+2}$ for
$ j\le\theta_i$. By the assumption, we have
$ a_{i,j}=a_j$, $b_{i,j}=b_j$, $\xi_{i,j}=\xi_j$ for $j\le\theta_i$.
Put $\theta_0=\min(\theta_1,\dots,\theta_s)$. Then $C(\tau)$ has the same
toric tangential direction of depth $\theta_0$ with $C(0)$ and we can
write $C^{(j)}(\tau)$ as:
$$\cases
C^{(j)}(\tau)&=\{(u_j,v_j)\in W_j; f_{\tau}^{(j)}(u_j,v_j)=0\}\\
f_{\tau}^{(j)}(u_j,v_j)&=
 (v_j^{a_{j+1}'}+\xi_{j+1}'u_j^{b_{j+1}'})^{A_{j+2}'}+
\text{(higher terms)}, \quad j\le \theta_0
\endcases\tag6.2.5
$$
where $A_{i+1}':=a_{i+1}'\cdots a_{\theta+1}'A_{\theta+2}'$ and
by the assumption, we have
$a_i'=a_i$, $ b_i'=b_i$,
$\xi_i'=\xi_i$, $A_{i+1}'=A_{i+1}$ for $i\le \theta_0$.
Comparing the defining equations of $C^{(i)}(\tau;1),\dots, C^{(i)}(\tau;s)$
and
$C^{i)}(\tau)$, we must have
$$
f^{(i)}_\tau(x,y)=\prod_{i=1}^s f_{\tau,i}^{(i)}(x,y)
,\quad
A_{1,i}+\cdots+ A_{s,i}=A_i,\quad i\le \theta_0
$$
As $A_{i,1}=a_{i,1}\cdots a_{i,\theta_i}A_{i,\theta_i+1}$
and $s\ge 2$, this implies that
$$\theta_i\le k-1\quad \text{and}\quad A_{i,\theta_i+1}<A_{\theta_i+1}.
\tag 6.2.6
$$
We use the following notations for simplicity.
$ \bar A_{j,i}:= A_{j,i}$ for $ i\le\theta_j+1$ and
$ \bar A_{j,i}:=0$ for $ i>\theta_j+1$.
Then by (6.2.5) and (6.2.6) we get
$$\sum_{j=1}^s \bar A_{j,i}
=A_i,~ i\le \theta_0+1\quad \text{and}\quad
\sum_{j=1}^s\bar A_{j,i}<A_i,~ i>\theta_0+1
\tag 6.2.7
$$
By Lemma (3.4.2), we have with
$P_{j,\theta_j+1}:=(a_{j,\theta_j+1},b_{j,\theta_j+1})$  that
$$\align
I(C(\tau;j),C(0);O)&=\sum_{i=1}^{\theta_j}a_ib_i A_{i+1}A_{j,i+1}
+I(P_{\theta_j+1},P_{j,\theta_j+1})
A_{\theta_j+2}A_{j,\theta_j+2}\\
&\qquad \quad \le \sum_{i=1}^{\theta_j}a_ib_i A_{i+1}A_{j,i+1}+
b_{\theta_j+1}A_{\theta_j+2}A_{j,\theta_j+1}
=\sum_{i=1}^{k}b_i A_{i+1}\bar A_{j,i}
\endalign
$$
Adding these inequalities for $j=1,\dots,s$ and using (6.2.6), we get
$$
I(C(\tau),C(0);O)\le \sum_{i=1}^{k}b_i A_{i+1}
\sum_{j=1}^s\bar A_{j,i}
\le \sum_{i=1}^{k}b_i A_{i+1}A_i,
$$
where the right side is  equal to $a_k I(C_{k-1},C(0);O)$ by
Theorem (4.5). Combining wit the assumption (6.1), we get
$I(C(\tau),C(0);O)=\sum_{i=1}^{k}b_i A_{i+1} A_i$.
This is equivalent to the following two equalities:
$$\align
&I(C(\tau;j),C(0);O)=\sum_{i=1}^{k}b_i A_{i+1}\bar A_{j,i},
\quad j=1,\dots,s\tag 6.2.8\\
&\sum_{j=1}^s\bar A_{j,i}=A_i,\quad i=1,\dots, k\tag6.2.9
\endalign
$$
By (6.2.7) and (6.2.6), (6.2.9) is equivalent to $\theta_0=k-1$.
Therefore (6.2.8) and (6.2.9) holds
 if only if
$\theta_i=k-1$ for $i=1,\dots,s$
and
$a_{k}b_{j,k}\ge a_{j,k}b_{k}$.
Assume that $a_{k}b_{j_0,k}>a_{j_0,k}b_{k}$ for some $j_0$.
Then we have:
$$\align
a_k I(C_{k-1},C(\tau);O)= \sum_{i=1}^{k-1}b_i A_{i+1}A_i
+a_k \sum_{j=1}^s b_{j,k}&A_{j,k+1}
>\sum_{i=1}^{k-1}b_i A_{i+1}A_i+
\sum_{j=1}^s{a_{j,k} b_k}A_{j,k+1}\\
& =\sum_{i=1}^{k}b_i A_{i+1}A_i=I(C(0),C(\tau);O)
\endalign
$$
Thus we get a contradiction $a_k I(C_{k-1},C(\tau);O)>I(C(0),C(\tau);O)$.
Therefore we must have
$a_{k}b_{j,k}=a_{j,k}b_{k}$, $ j=1,\dots,s$.
As $\gcd(a_k,b_k),~\gcd(a_{j,k},b_{j,k})=1$, this is possible if and only if
$a_{j,k}=a_k$ and $b_{j,k}=a_k$. Again this gives a contradiction:
$A_k=A_{1,k}+\cdots+A_{s,k}=sA_k.$
This proves  the irreducibility of $C(\tau)$ and
the proof of Theorem (6.2) is now completed. ~\QED
\enddemo
\ninpar
\S 7. {\bf An example of an equi-singular family}.
We study a typical equi-singular family $f_t(x,y):=f(x,y)+tx^m,$
where $f(x,y)$ is a polynomial, whose Newton diagram
$\De(f;(x,y))$
is a triangle with the vertices $A=(0,n),B=(b_1A_2,0),C=(m,0)$ with $m>b_1A_2,$
 having the initial expansion
$f(x,y)=(y^{a_1}+\xi_1 x^{b_1})^{A_2}+
\text{(higher terms)},\quad a_1\ge 2,$
and defining an irreducible germ of a plane curve
$C=\{ f(x,y)=0\}$
 at the origin.
Then the $a$-th \TS~ polynomial of $f_t(x,y)$
does not depend on $t$ for any $a|n$ with $1 < a$,
 so we can apply the previous consideration to the family
of germes
$
C(t):=\{(x,y)\in \bC^2; f_t(x,y)=0\}$.
A similar family is studied by Ephreim [6] using polar invariants. Let
$\{P_i={}^t(a_i,b_i);i=1,\dots, k\}$
be the weight vectors of the  \Ts ~resolution tower. Let $h_i$
be the  $A_{i+1}$-th \TS~ polynomials of $f_t(x,y)$ for $i=1,\dots, k-1$
and let $C_{i}=\{(x,y)\in \bC^2; h_i(x,y)\}$.
\proclaim{Proposition (7.1)} With the above assumptions and notations, we have
$I(C(t),C(s);O)=nm$ for $ t\ne s$ and
$a_k I(C_{k-1},C(t);O)\le nm$ for any $t\in \bC$.
\endproclaim
\demo{Proof}
For the proof of the equality, note:
$I(C(t),C(s);O)=\dim_{\bC}\bC \{x,y\}/(f_s, f_t)$ and therefore
$=\dim_{\bC}\bC \{x,y\}/(f_t, (t-s)x^m)=nm$.
To prove the inequality, we first observe the Newton diagram
$\De(h_{k-1})$ is a subset of the triangle $\De'$
whose vertices are $A'=(0,n/a_k),B'=(b_1A_2/a_k,0),C'=(m/a_k,0)$.
In the case of $m=n$, the assertion follows from  the Bezout theorem
in $\bP^2$:
$a_k I(C_{k-1},C(t);O)\le a_k \bar C_{k-1}\cdot \bar C(t)=n^2=nm$,
where $\bar C$ is the projective compactification of $C\subset \bC^2$
and the right side is the  intersection number in $\bP^2$.
In the case $m\ne n$, we need another argument.
Choose a small ball $B$ centered at the origin
$B$   containing no other intersection than the origin $O$.
Let $f_s'(x,y)=f_s(x,y)+\eps_1$ and $h_{k-1}'(x,y)=h_{k-1}(x,y)+\eps_2$
and let
$C(s)'=\{(x,y)\in \bC^2;f_s'(x,y)=0\}$ and
$C_{k-1}'=\{(x,y)\in \bC^2;h_{k-1}'(x,y)=0\}$.
For sufficiently small $\eps_1,\eps_2$
 the intersection $C(s)'\cap C_{k-1}'$ is a subset of the torus
$\torus{2}$ and the number  of the points counted with  multiplicity of
$C(s)'\cap C_{k-1}'$
in $B$ is equal to $I(C_{k-1},C(t);O)$.
 The Newton diagram $\De:=\De(f_s';(x,y))$ is the triangle with vertices $O,A$
and $C$.
The number of intersection points $C(s)'\cdot C_{k-1}'$ in $\bC^{*2}$
 is bounded by
 the theorem of Bernshtein ([5],[28]):
 $$
C(s)'\cdot C_{k-1}'=2V_2(\De(f_s'),\De(h_{k-1}'))
\le 2V_2((\De,\De/a_k)=2\Vol(\De)/a_k=nm/a_k
$$
Here $V_2(\De_1,\De_2)$ is the  Minkowski's mixed volume
and we have used the monotone increasing property
 of the  Minkowski's mixed volume  to
the inclusion $\De(h_{k-1}')\subset \De/a_k$.
See [5, 28]. As  $C(s)'\cdot C_{k-1}'\ge I(C_{k-1},C(s))$,
the inequality of the proposition follows. ~\QED
\enddemo
\ninpar
\S 8. {\bf The equi-singularity at infinity and
the Abhyankar-Moh-Suzuki theorem.}
\nl
 Let
 $F:\bC^2\to \bC$ be a polynomial mapping of degree $n$.
We say that $\tau\in\bC$ is
a {\it regular value at infinity} if there exits a large number
$R$ and a positive number $\de$ so that the restriction
$F:E_\infty(R,\de)\to D_\de$ is a trivial fibration
where
$$D_\de=\{\eta\in \bC; |\eta-\tau|\le \de\},\qquad
E_\infty( R,\de)=\{(x,y);F(x,y)\in D_\de,~\sqrt{|x|^2+|y|^2} >R\}$$
Let $C_t=F\inv(t)$ and let $\bar C_t$ be the projective
compactification. The set
$\bar C_t  -  C_t=\{\rho_1,\dots, \rho_\ell\} \subset L_\infty$ be
 does not depend on $t$. We recall that
\proclaim{Proposition (8.1)([10])}
$\tau\in \bC$ is a regular value at infinity  if and only if
the family of germs of plane curves $\{(\bar C_t,\rho_i);t\in \bC\}$ is
topologically stable at $t=\tau$ for any $i=1,\dots, \ell$.
\endproclaim
We consider hereafter the  simplest case that $C_0$ has one place at infinity,
 say at $\rho=(1;0;0)$. Namely we assume that
$\ell=1$ and the germ $(\bar C_0,\rho)$ is irreducible.
Then
$F(x,y)$ is written as
$$F(x,y)=(y^{a_1}+\xi_1 x^{c_1})^{A_2}+
\text{(lower terms)},\quad c_1< a_1, n=a_1A_2.\tag 8.1.1$$
for some integers $a_1\ge 2$, $A_2$ and  $1\le c_1 < a_1$.
As $\bar C_0\cap L_\infty=\{\rho\}$ and $\bar C_0$
is assumed to be  locally irreducible at
$\rho_1$, the polynomial $F(x,y)$ has only one outside  face
and its outside face function has
only one factor. See [18] or [26].
The standard affine coordinates
$ u=Z/X,~v=Y/X$
are centered at $\rho$ and
the curve $\bar C_t$ is defined by
$\{f_t(u,v)=0\}$ where $f_t(u,v)= f(u,v)-t u^n$
and $f(u,v)=F(1/u,v/u)\times u^n$.
In this simplest case, we have the initial expansion
$$
f(u,v)=(v^{a_1}+\xi_1 u^{b_1})^{A_2}+\text{(higher terms)}\quad
\tag 8.1.2
$$
where $b_1=a_1-c_1$.
Let $C_t^\infty=\{(u,v)\in \bC^2;f(u,v,t):=f(u,v,t)=f(u,v)-t u^n=0\}$.
We can apply Theorem (6.2) using Proposition (7.3) to this family and
we obtain:
\proclaim{Theorem (8.2)} For the mapping $F$ and the family
$\{(C_t^\infty,O);t\in \bC\}$ holds the following
:\nl
{\rm (1)} The family of germs of plane curves $\{(C_t^\infty,O);t\in \bC\}$
is an equi-singular family of irreducible curves and the \TS~
resolution tower of $(C_0^\infty,O)$ resolves simultaneously
each curve of the family $\{(C_t^\infty,O);t\in \bC\}$.
\nl
{\rm (2)} The mapping $F:\bC^2\to \bC$ has no critical point at
infinity.
\endproclaim Ephraim has also obtained  a similar result
about the equi-singularity using a different method ([6]).
See also Moh [20].
\par
Before giving  applications, we will need the
following facts. Let $D\subset \bP^2$ be a projective curve of degree $n$
and let $q_1,\dots, q_\nu$ be the singular points of $D$. Then
by Pl\"ucker's formula and by Mayer-Vietoris argument,
the topological Euler number
of $D- \{q_1,\dots, q_\nu\}$ is given by
$$\chi(D - \{q_1,\dots, q_\nu\})=2-\nu-(n-1)(n-2)+
\sum_{i=1}^\nu \mu(D;q_i).$$
{}From this equality follow two equivalences.
First, $\mu(D;q_1)=(n-1)(n-2)$ if and only if the curve $D-\{q_1\}$
is smooth  and homeomorphic to the line $\bC$. Second,
$\mu(D;q_1)=(n-1)(n-2)-2g $ and  $\nu=1$
 if and only if the curve $D-\{q_1\}$
is smooth  and homeomorphic to a punctured Riemann surface of genus $g$.
As a first application, we will give an elementary proof
of:
\proclaim{Theorem (8.3) (Abhyankar-Moh [4], Suzuki [29])}
Let $F(x,y)$ be a polynomial of two variables
of degree $n$ and assume that the plane curve
$C=\{(x,y)\in \bC^2;F(x,y)=0\}$ is smooth and homeomorphic to
the complex line $\bC$. Then there exists another polynomial
$G(x,y)$ so that $(F,G)$ is an automorphism of $\bC^2$.
\endproclaim
\demo{Proof}
 The polynomial $F(x,y)$ has one place at infinity, say at $\rho=(1;0;0)$.
To prove  the theorem  by the induction on $n=\degree F(x,y)$,
it is enough to show that
$c_1=1$
in (8.1.1). In fact, if $c_1=1$,
 we apply the coordinate
change
$(X,Y)=(y^{a_1}+\xi_1 x,y)$ and achieve
$\deg F((X-Y^{a_1})/\xi_1,Y)<n$. Therefore
the assertion is proved by the induction on $\deg F$.
\par
Let $C_0^\infty$ and $\bar C_0$ be as above.
We have
$\mu( C_0^\infty;O)=(n-1)(n-2)$ since the smooth part of the
curve $\bar C_0$ is homeomorphic to the line $\bC$.
Let us consider the \TS~ resolution tower of $(C_0^\infty,O)$:
$\Cal T=\{X_k\mapright{p_k}X_{k-1}\to \cdots\to X_1\mapright{p_1}X_0=\bC^2\}$
 and let
$P_i={}^t(a_i,b_i)$, $i=1,\dots, k$
be the weight vectors of the tower. Then by Theorem (5.1) and $n=A_1$, we have
$$(A_1-1)(A_1-2)=\mu(C_0^\infty,O)=1-A_1+\sum_{i=1}^k(A_{i}-1) b_iA_{i+1},$$
which leads to
$$
\sum_{i=1}^k(A_{i}-1) b_iA_{i+1}=(A_1-1)^2.\tag a
$$
{}From Theorem (4.5) and  Bezout theorem, we deduce
$$
\sum_{i=1}^k a_ib_i A_{i+1}^2\le A_1^2\tag b
$$
since
$a_k I(C_{k-1},C(0);\xi_0)=\sum_{i=1}^k a_ib_i A_{i+1}^2
\le  a_k\bar C_{k-1}\cdot \bar C(0)=A_1^2$.
Now we are ready to show $c_1=1$. We follow the proof of
Abhyankar-Moh, Lemma (3.1), [4].
Recall that $c_1=a_1-b_1$.
For the case $k=1$, the equality (a) reads
$(a_1-1)b_1=(a_1-1)^2$. Thus we get $c_1=a_1-b_1=1$. For the case $k\ge 2$,
we rewrite (a) and (b)  as
$$\align
\sum_{i=2}^k(A_{i}-1) b_iA_{i+1}&=(a_1A_2-1)(c_1A_2-1)\tag c\\
\sum_{i=2}^k a_ib_i A_{i+1}^2&\le c_1a_1A_2^2\tag d\endalign
$$
Thus taking the sum : (c)$\times A_2~+~$(d)$\times (1-A_2)$,
we obtain
$$
\sum_{i=2}^k b_iA_{i+1}(A_i-A_2)\ge A_2^2((a_1-1)(c_1-1)-1)+A_2.$$
The left side is obviously negative. The right hand is negative
 only if $c_1=1$, which completes the proof.~\QED
\enddemo
\proclaim{Theorem (8.4)} The weight vectors of a good toric
resolution of the singularity at infinity of a smooth
acyclic curve in $\bC^2$ satisfy
 $b_i=a_{i-1}a_i-1$ for  each $i=1,\dots,k$
where $a_0=1$.\endproclaim
\demo{Proof}
Substituting $c_1=1$ in (c) and (d), we get
$$\align
\sum_{i=3}^k(A_{i}-1) b_iA_{i+1}&=(a_2A_3-1)(c_2A_3-1)\tag e\\
\sum_{i=3}^k a_ib_i A_{i+1}^2&\le c_2a_2A_3^2\tag f\endalign
$$
where $c_i=a_{i-1}a_i-b_2$ for $i=2,\dots,k$.
Thus again taking the sum :
(e)$\times A_3~+~$(f)$\times (1-A_3)$,
we obtain
$$
\sum_{i=3}^k b_iA_{i+1}(A_i-A_3)\ge A_3\{((a_2-1)(c_2-1)-1)A_3+1\}$$
The left side is obviously negative. The right hand is negative
 only if $c_2=1$. The assertion for $i\ge 2$
 can be proved by an easy induction.~\QED
\enddemo
\par
The following example shows that all weight vectors having the
property of  Theorem (8.4)  occur.
\example{Example (8.5)}
Let $a_i\ge 2,~i=1,\dots,k$ be given integers, and let $n=a_1\cdots a_k$.
Let us consider the sequence of automorphisms:
$$
\varphi_i:
\pmatrix x_{i}\\ x_{i+1}\endpmatrix \mapsto\pmatrix x_{i+1}\\
x_{i+2}\endpmatrix,\quad x_{i+2}=x_i+x_{i+1}^{a_{i+1}},\quad
i=0,\dots, k-1
$$
where $x_0=x$ and $x_1=y$.
Let $F(x,y)=x_{k+1}(x,y)$. Then $F(x,y)$ obviously
 satisfies the assumption of Theorem (8.3).
Let
$$\cases
h_0&=v\\
h_1(u,v)&=v^{a_1}+u^{a_1-1}\\
h_2(u,v)&=h_1(u,v)^{a_2}+h_0(u,v)u^{a_1a_2-1}\\
h_i(u,v)&=h_{i-1}^{a_{i}}+h_{i-2}(u,v)u^{a_1\cdots a_{i-2}(a_{i-1}a_i-1)},
\quad 2\le i\le  k
\endcases\tag $\sharp_1$
$$
and let $f(u,v)=h_k(u,v)$.  Then $(C_0^{\infty},O)$
is defined by
$C_0^{\infty}=\{(u,v)\in\bC^2; f(u,v)=0\}$.
It is easy to see that $h_i$ is the $A_{i+1}$-th
\TS~ polynomial of $f$. By an inductive argument we can prove
that the weight
vectors of the \TS~ resolution tower is given by
$$P_1=\pmatrix a_1\\a_1-1\endpmatrix,\quad
 P_2=\pmatrix a_2\\a_1a_2-1\endpmatrix,\dots,
P_k=\pmatrix a_k\\a_{k-1}a_k-1\endpmatrix,$$
and   the pull-backs of the \TS~ polynomials to $X_i$ are given by
$$\cases
\Phi_i^* h_i(u_i,v_i)&=u_i^{m_i(h_i)}v_i\\
\Phi_i^* h_{i+1}(u_i,v_i)&=u_i^{m_i(h_{i+1})}(v_i^{a_{i+1}}+\xi_{i+1}
u_i^{a_{i}a_{i+1}-1})\\
\Phi_i^* h_{j}(u_i,v_i)&=u_i^{m_i(h_{j})}(\bar h_{j-1}^{a_j}+\xi_{i+1}
\bar h_{j-2}u_i^{a_i\cdots a_{j-2}(a_{j-1}a_j-1)})
\endcases\tag $\sharp_{i+1}$
$$
where
$\bar h_j(u_i,v_i):=\Phi_i^* h_{j}(u_i,v_i)/u_i^{m_i(h_{j})}$
and $\xi_{i+1}$ is a unit in a neighbourhood $W_i$  of $\Xi_i$.
The Milnor number is equal to $(n-1)(n-2)$ by Theorem (5.1).
 It is convenient to introduce the notation $a_0=1$ and
$m_0(h_i)=0$ to understand ($\sharp_1$)
as a special cases of ($\sharp_{i+1}$).
\endexample
\remark{ Remark (8.6)} Let $F(x,y)$ be a polynomial of degree $n$ and
coefficients in a subfield $k$ of $\bC$, such that the curve
$C=\{F(x,y)=0\}\subset \bC^2$
is smooth and contractible. Then the completion of $C$ requires one extra
 point $\rho$ at infinity  having its coordinates in $k$. So, after a
linear change of coordinates defined over $k$,  the pencil
$L_t=\{ y=t\}_{t\in \bC}$
passes through $\rho$. Let $(B(t),t)$ be   the barycenter, computed
in the affine line $L_t$, of  the points of the intersection
$L_t\cap C$, weighted by the  multiplicity.  The automorphism
$(x,y) \to (x-B(y),y)$, which is
 defined over $k$, moves the curve $C$ to a curve $C'$ of lower degree
and having at infinity one Puiseux pair less.
In the notation of Theorem (8.3),
we can write $B(y)=-y^{a_1}/\xi_1+\text{(lower terms)}$.
 The iteration of this
procedure constructs an automorphism defined over $k$, which moves
the curve $C$ to a line.
Of course, we can apply this procedure to any curve $D=\{G(x,y)=0\}$,
as long as  the completion of the curve $D$ has only one
irreducible singularity at infinity and $c_1=1$.
After at most $\log_2(\degree(G))$
automorphism applied to  the curve $D$,
either the curve $D$ becomes a line and the equation linear,
in which case the curve $D$ was smooth and contractible,
or the curve $D$  becomes a curve for which $c_1\ge 2$.
This provides a  test for
the contractibility and smoothness of the curve $D$.
\endremark
We can apply the above  remark and  argument to get:
\proclaim{Theorem (8.7)} Let $C \subset \bC^2$ be a smooth curve
homeomorphic to a  surface with one puncture of genus
$g$,~$g=1$ or $2$.
Then there exist an automorphism of $\bC^2$ moving the curve $C$ to
a  smooth  cubic curve  which is tangent to the line at infinity
with the intersection multiplicity 3 if $g=1$,  and
to  a curve of degree 5 with a non-degenerate
singularity at infinity, if $g=2$.
\endproclaim
\demo{Proof}Let
$\Cal T=\{
X_k\mapright{p_k} X_{k-1}\to\cdots \to X_1\mapright{p_1} X_0=\bC^2\}$
be a Tschirnhausen tower of resolution of the singularity at infinity
with the corresponding weight vectors
$\{P_i={}^t(a_i,b_i);i=1,\dots, k\}$  as before.
 Applying barycentric automorphisms  if necessary,
we can assume that $c_1\ge 2$ (see Remark (8.6)).
The equalities (a) and (c) take the following form.
$$
\align
\sum_{i=1}^k(A_{i}-1) b_iA_{i+1}&=(A_1-1)^2 -2g\tag $a_g$\\
\sum_{i=2}^k(A_{i}-1) b_iA_{i+1}&=(a_1A_2-1)(c_1A_2-1)-2g\tag $c_g$
\endalign
$$
The inequality (b) and (d) are valid as before.
Then taking the sum : $(c_g)\times A_2~+~$(d)$\times (1-A_2)$,
we obtain
$$0 \ge \sum_{i=2}^k b_iA_{i+1}(A_i-A_2)\ge A_2^2((a_1-1)(c_1-1)-1)
- (2g-1)A_2. \tag $e_g$
$$
Note also that  $a_1\ge 3$ by the assumption $c_1\ge 2$.
\nl
(1) Assume first $k=1$. Then we have
$(a_1-1)(c_1-1)=2g$. So, for the natural number $c_1:=a_1-b_1$ we have
$a_1 > c_1=1+2g/(a_1-1)$. We conclude that $a_1=3$, if $g=1$ and that $a_1=5$,
if $g=2$.  \nl
So, for $g=1$, we have $a_1=3$,
$c_1=2$ . As $b_1=1$,
the curve $C$ has no singularity at infinity but $C$ is tangent
to the line at infinity with the tangent multiplicity
3. An example of such curve is given by
$C=\{y^3+x^2+1=0\}$.
For $g=2$ we have $a_1=5$ and $c_1=2$.
An example of such curve is given by $C=\{y^5+x^2+1=0\}$.
The curve $C$ has a non-degenerate cusp singularity at infinity.
\nl
(2) Now we show the case $k\ge 2$ does not occur. \nl
With $a_1 > c_1 \ge 2$,  we  deduce  from the
inequalities $(e_g)$:
$$A_2\le \frac{2g-1}{((a_1-1)(c_1-1)-1)} \le  (2g-1) \leqno (\star)$$
If $g=1$,  we get, from $(\star)$, $A_2=1$, and hence $k=1$.
\nl
If $g=2$, we reduce from $(\star)$: $ A_2\le 3$
that  $A_2=1,2$ or $3$.
We first rule out the case $A_2=3$: indeed, from $(\star)$
 we conclude
$k=2,~a_1=3,~c_1=2,~b_1=1$. So, $a_2=A_2=3,~n=9$ and $b_2=18$ by $(a_g)$. This
is not possible since  we have assumed
$\gcd(a_2,b_2)=1$.\nl
Next, we rule out the case  $A_2=2$: from $(\star)$,
we conclude $k=2,~a_1=3,~c_1=2,~b_1=1$. So, $a_2=A_2=2,~n=6,~b_1=1$ and
$b_2=11$ by $(a_g)$.
Thus the tower has the weight vectors
$P_1={}^t(3,1)$ and $P_2={}^t(2,11)$. No easy contradiction yet.
However
we assert that there is no polynomial
$f(u,v)$ of degree 6, irreducible at the origin, whose weight vectors are
as above. Indeed, let
$f(u,v)=(v^3+ u)^2+\sum_\nu c_\nu u^{\nu_1}v^{\nu_2}$
where
$6<3\nu_1+\nu_2,\quad \nu_1+\nu_2\le 6$.
Consider an admissible toric modification $p:X_1\to \bC^2$. We may assume that
$\si=\Cone(E_1,P_1)$ is the left toric  cone of the divisor $\what E(P_1)$
 and let $(s,t)$ be the toric coordinates. Then we have
$u=st^3$ and $v=t$. The pull backs can be written as
$\pi_\si^*(v^3+ u)^2= t^6(1+s)^2$ and
$\pi_\si^* u^{\nu_1}v^{\nu_2}=s^{\nu_1}t^{3\nu_1+\nu_2}$.
So, in $t^{-6}\pi_\si^*f(s,t)$ the monomial $t^{11}$ does not
occur and hence $P_2$ is not the second weight vector for $f(u,v)$.
Thus  this case does not occur. So, $A_2=1$ proving $k=1$.~\QED
\enddemo
\remark{Remark (8.9)}
Using $(\star)$ and inequality:  $2^{(k-1)}\le A_2$,
we  get the following estimate for
 the length of the tower:
$$k\le \log_2(2g-1)+1,\qquad \text{for}\quad g\ge 1,~c_1\ge 2.$$
The classification for $g\ge 3$ is more complicated, as
the model is not unique. For example, in the case
of $g=3$, we can move $C$ by an  automorphism to one of the following.
\nl
(a) $k=1,~P_1={}^t(4,1)$ and $n=4$. The curve is smooth at infinity
and tangent to the line at infinity at a single point.
An example is given by
$y^4+x^3+1=0$.\nl
(b) $k=1,~P_1={}^t(7,5)$ and  $n=7$. The curve has a non-degenerate cusp
 singularity
at infinity. An example is given by
$y^7+x^2+ 1=0$.
\nl
(c) $k=2,~P_1={}^t(3,1),~P_2={}^t(2,9)$ and $n=6$. An example is given by
$(y^3+x^2)^2+ x$.
\endremark
Professor M. Miyanishi recently communicated to us
that he gave a new proof of Theorem (8.3)
using the classification of surfaces [8].
Also, the paper [32] contains interesting results about contractible affine
curves with one isolated singularity.

\enddocument